\documentclass[a4paper,11 pt,reqno]{amsart}
\usepackage{a4,amsbsy,amscd,amssymb,amstext,amsmath,latexsym,oldgerm,enumerate,verbatim,graphicx,xy}

\makeatletter                                       
\def\l@subsection{\@tocline{2}{0pt}{2.5pc}{2.5pc}{}}
\makeatother                                        

\makeatletter                                                  
\def\chapter{\clearpage\thispagestyle{plain}\global\@topnum\z@ 
\@afterindenttrue \secdef\@chapter\@schapter}
\makeatother

\newtheorem{thm} {Theorem} [section]
\newtheorem{prop}{Proposition} [section]
\newtheorem{lem} {Lemma} [section]

\newtheorem{corgl} {Corollary}

\newtheorem{cornn}{Corollary}

\theoremstyle{definition}

\newtheorem{rem} {Remark} [section]
\newtheorem{rems} [rem]{Remarks}

\newcommand{\mf}{\mathfrak}
\newcommand{\mc}{\mathcal}
\newcommand{\mb}{\mathbb}

\newcommand{\ov}{\overline}
\newcommand{\un}{\underline}

\newcommand{\sm}{\setminus}         

\newcommand{\ot}{\otimes}           
\newcommand{\la}{\langle}
\newcommand{\ra}{\rangle}

\newcommand{\End}{{\rm End}}

\newcommand{\Mat}{{\rm Mat}}

\newcommand{\ind}{{\rm ind}}

\newcommand{\Spec}{{\rm Spec}}

\newcommand{\Sym}{{\rm Sym}} 

\newcommand{\rk}{{\rm rk}}

\let\ttie\t
\newcommand{\tie}[1]{{\let\t\ttie \ttie#1}}
\renewcommand{\t}{\mf{t}}  

\newcommand{\GL}{{\rm GL}}

\newcommand{\Ort}{{\rm O}}

\newcommand{\Sp}{{\rm Sp}} 

\newcommand{\stimes}{\mbox{\fontsize{7}{6}\selectfont $\times $}}
\newcommand{\e}{\epsilon}
\newcommand{\ve}{\varepsilon}

\newcommand{\SpM}{{\rm SpM}}
\newcommand{\GSp}{{\rm GSp}}
\newtheorem{mycor}{\addtocounter{corgl}{1}Corollary~\thecorgl\ to Proposition~\ref{prop.extbasis}}

\xyoption{all}

\begin{document}

\title{A bideterminant basis for a reductive monoid}

\begin{abstract}
We use the rational tableaux introduced by Stembridge to give a bideterminant basis for a normal reductive monoid and for its variety of noninvertible elements. We also obtain a bideterminant basis for the full coordinate ring of the general linear group and for all its truncations with respect to saturated sets. Finally, we deduce an alternative proof of the double centraliser theorem for the rational Schur algebra and the walled Brauer algebra over an arbitrary infinite base field which was first obtained by Dipper, Doty and Stoll.
\end{abstract}

\author[R.\ H.\ Tange]{Rudolf Tange}


\maketitle

\section*{Introduction}
Let $k$ be an infinite field. Assume for the moment that $k=\mb C$. Amongst the several instances of Schur-Weyl duality there are:
\begin{enumerate}[$\bullet$]
\item The symmetric group $\Sym_r$ and the general linear group $\GL_n$ acting on $V^{\otimes r}$, $V=k^n$ the natural module of $\GL_n$.
\item The Brauer algebra $B_r(n)$ or $B_r(-n)$ and the orthogonal or symplectic group acting on $V^{\ot r}$, $V=k^n$ the natural module of the orthogonal or symplectic group.
\item The walled Brauer algebra $B_{r,s}(m)$, see \cite{BCHLLS} or \cite{CdVMD}, and the general linear group $\GL_m$ acting on $V^{\ot r}\ot(V^*)^{\ot s}$, $V=k^m$ the natural module of $\GL_n$.
\end{enumerate}
The initial motivation of this paper was the third instance of Schur-Weyl duality for $k$ any infinite field. The way to understand this duality is to consider the action of the  orthogonal group $\Ort_m\subseteq\GL_m$ as well. For the orthogonal group we have $V=V^*$. So, by the second instance, we should have that the centraliser algebra is a subalgebra of the Brauer algebra $B_{r+s}(m)$. The problem is to show that the image of the walled Brauer algebra and the enveloping algebra of $\GL_m$ in $\End_k(V^{\ot r}\ot(V^*)^{\ot s})$ are each other's centraliser. The hard part here is to show that the centraliser algebra of the walled Brauer algebra is the enveloping algebra of $\GL_m$. This led us to study a certain reductive submonoid $M$ of $\Mat_m\times\Mat_m$ which is the main topic of this paper. In fact we have to study $M$ as a monoid scheme and then deduce afterwards that it is reduced over $k$. To do this we will apply modified versions of the methods of \cite{T} and \cite{Oe}. We will state results in a form which makes the link with the general theory of reductive groups clear. This is made possible by an improved result about straightening in \cite{DeCEP}. To make the exposition as clear as possible we first consider bideterminant bases for the variety of $n\times n$-matrices $\Mat_n$ and for the symplectic group $\Sp_{2m}$. Our method is to move in the sequence $\Mat_n, \Sp_{2m},M$ things as much as possible in the $\Mat_n$-direction and to give for the other cases only proofs if they differ significantly from the previous case.

The paper is organised as follows. In Section~\ref{s.prelim} we state the results on bideterminants and straightening from \cite{DeCEP} and mention the link with the general theory of reductive groups. Furthermore we give a result which relates enveloping algebra of a submonoid of $\Mat_m\times\Mat_m$ in $\End_k(V^{\ot r}\ot V^{\ot s})$ with its vanishing ideal. In Section~\ref{s.symplecticstraightening} we give the basic results about bideterminants and straightening in the symplectic case. This is based on Oehms' work \cite{Oe}. Since we want to explain the link with the general theory of reductive groups, we have to do the straightening directly in $k[\SpM_n]$ and not in $k[\SpM_n]/(d)$ as in \cite{Oe}. The point is that we have to know something about the shapes of bideterminants of lower degree that show up during straightening. In Section~\ref{s.bideterminant basis} we finally give the results about bideterminants and straightening for the monoid $M$. To construct our bideterminant basis we need rational tableaux as introduced by Stembridge \cite{St}. We also show that the full coordinate ring of $\GL_m$ has a bideterminant basis and we show that $M$ is normal and therefore, by a result of Rittatore, Cohen-Macaulay. In Section~\ref{s.doublecentraliser} we give a proof of the double centraliser theorem for $\GL_m$ and the walled Brauer algebra $B_{r,s}(m)$ acting on $V^{\ot r}\ot(V^*)^{\ot s}$.

\section{Preliminaries}\label{s.prelim}
Throughout this paper $k$ denotes an infinite field and $K$ denotes the algebraic closure of $k$. We will denote a scheme over a commutative ring $R$ always like $X_R$, the base ring as a subscript. All schemes in this paper are affine and algebraic over the base ring. We refer the reader to \cite{Jan} or \cite{DemGa} for the basics about schemes. We will only make a very modest use of these. If $A$ is a commutative $R$-algebra, then we write $X(A)=X_R(A)$ for the points of $X_R$ over $A$. {\it In case $R$ is $\mb Z$ or our infinite field $k$ we will, to keep notation manageable, simply denote $X(k)$ by $X$.} We denote the coordinate ring of $X_R$ by $R[X_R]$, it can be identified with the $R$-algebra of morphisms $X_R\to\mb A^1_R$, where $\mb A^1_R$ denotes the affine line over $R$. There is a canonical homomorphism of $R$-algebras from $R[X_R]$ to the algebra of $R$-valued functions on $X(R)$; we denote its image by $R[X(R)]$. If the base ring is $\mb Z$ or $k$, then the epimorphism $k[X_k]\to k[X]$ is an isomorphism if and only if $X_k$ is reduced and $X(k)$ is dense in $X(K)$. To avoid artificial generality, we will work over $\mb Z$ if we want results valid for more general rings than fields. The reader can obtain the result he is interested in by tensoring with his favourite ring (e.g. $\mb C$).

Now let $n$ be an integer $\ge 1$, let $\Mat_{n,\mb Z}$ be the $\mb Z$-scheme of $n\times n$-matrices and let $\GL_{n,\mb Z}$ be the $\mb Z$-group scheme of $n\times n$-matrices. So, for a commutative ring $A$, $\Mat_n(A)$ is the set of $n\times n$ matrices over $A$ and $\GL_n(A)$ is the group of invertible $n\times n$ matrices over $A$. We have $\mb Z[\Mat_{\mb Z}]=\mb Z[(x_{ij})_{ij}]$, the polynomial algebra over $\mb Z$ in the matrix entries $x_{ij}$. Furthermore, $\mb Z[\GL_{\mb Z}]=\mb Z[\Mat_{\mb Z}][\det^{-1}]$. The group scheme $\GL_{n,\mb Z}\times \GL_{n,\mb Z}$ acts on $\Mat_{n,\mb Z}$ via $(g,h)\cdot A=gAh^{-1}$ for $g,h\in\GL_n(R)$, $A\in\Mat_n(R)$ and $R$ a commutative ring. We obtain an action on the coordinate ring of $\Mat_{n,\mb Z}$ which is given by $((g,h)\cdot f_R)(A)=f_R(g^{-1}Ah)$ for $f\in\mb Z[\Mat_{\mb Z}]$, $g,h\in\GL_n(R)$, $A\in\Mat_n(R)$ and $R$ a commutative ring. The action of the left resp. right factor of $\GL_{n,\mb Z}\times \GL_{n,\mb Z}$ on $Z[\Mat_{n,\mb Z}]$ comes from the left resp. right multiplication and therefore we refer to it as the {\it left} resp. {\it right regular action}. We note here that for $R$ any commutative ring $R[\GL_R]$ is flat over $R$, since $R[\Mat_{n,R}]$ is free over $R$ and localisation is exact. As a consequence the category of $\GL_{n,R}$-modules is abelian, see \cite[I.2.9]{Jan}. In fact one can show with a bit more effort that $R[\GL_{n,R}]$ is free over $R$. In Section~\ref{s.bideterminant basis} we will see an explicit basis.

We denote the character group (the homomorphisms to $\GL_{1,\mb Z}$) of the maximal torus of $\GL_{n,\mb Z}$ of diagonal matrices by $X$. We have $X\cong \mb Z^n$ where the $i^{\rm th}$ diagonal matrix entry corresponds to the $i^{\rm th}$ standard basis element $\ve_i$ of $\mb Z^n$.
We denote the set of dominant (relative to the Borel subgroup of upper triangular matrices) weights $\lambda_1\ge\lambda_2\ge\cdots\ge\lambda_n$ by $X^+$ and the set of the polynomial dominant weights, i.e. partitions of length $\le n$, $\lambda_1\ge\lambda_2\ge\cdots\ge\lambda_n\ge0$, by $\Lambda^+(n)$. Put $\un n:=\{1,\ldots,n\}$. For $\lambda\in X$ we put $|\lambda|=\sum_{i=1}^n\lambda_i$. In case $\lambda$ is a partition we say that $\lambda$ is a partition of $r$ if $|\lambda|=r$. Furthermore, $|\{i\in\un n\,|\,\lambda_i\ne0\}|$ is called the length of $\lambda$. For $\lambda,\mu\in X$ we write $\mu\le\lambda$ if $\lambda-\mu$ is a sum of positive roots, i.e. if $|\lambda|=|\mu|$ and $\sum_{i=1}^j\mu_i\le\sum_{i=1}^j\lambda_i$ for all $j\in\un n$. For partitions this is the well-known dominance order. Now let $\lambda$ be a partition of length $l\le n$. Then the shape or Young diagram of $\lambda$ is the set of pairs $(i,j)$, where $1\le j\le\lambda_i$ and $1\le i\le l$. The conjugate partition of $\lambda$ is the partition whose diagram is the transpose of that of $\lambda$. We denote it by $\lambda'$. Note that the length of $\lambda'$ is $\lambda_1$, the number of columns of $\lambda$.

We assume given a linear order $\preceq$ on $\un n$, for example the natural order. A {\it tableau of shape $\lambda$} is a function from the shape of $\lambda$ to $\un n$ and a {\it bitableau of shape $\lambda$} is a pair $(S,T)$ where $S$ and $T$ are tableaux of shape $\lambda$. A tableau is called {\it standard} if, according to $\preceq$, its entries are strictly increasing down the columns and weakly increasing in the rows from left to right. Note that for two linear orderings $\preceq_1$ and $\preceq_2$ of $\un n$ there is a permutation of $\un n$ which induces a bijection between the tableaux that are standard relative to $\preceq_1$ and those that are standard relative to $\preceq_2$. For a partition $\lambda$ of length $l\le n$ we define the {\it canonical tableau of shape $\lambda$}, $T_\lambda$, to be the tableau of shape $\lambda$ whose entries in the $i^{\rm th}$ row are all equal to $i$, $1\le i\le l$. Now let $(S,T)$ be a bitableau of shape $\lambda$. Then the product of the minors
\begin{equation}\label{eq.bideterminant}
\det\big( (x_{S(r,i)\,T(s,i)})_{1\le r,s\le\lambda_i'}\big)\, ,
\end{equation}
$1\le i\le\lambda_1$, in $\mb Z[\Mat_{n,\mb Z}]$ is called the {\it bideterminant} associated to $(S,T)$ and we denote it by $(S\,|\,T)$. So, as in \cite{Green}, we form bideterminants according to pairs of columns in a bitableau rather than to pairs of rows as in \cite{DeCEP}. Put differently, the bideterminant associated to $(S,T)$ in \cite{DeCEP} is the bideterminant that we associate to $(S',T')$, where $S'$ and $T'$ are the transposed tableaux of $S$ and $T$ (they have shape $\lambda'$). The reader should note that the bideterminants associated to bitableaux of shape $r\ve_1$ for some $r\ge0$ are precisely the monomials in the $x_{ij}$.
We define the {\it content} or {\it weight} of a tableau $T$ to be $\sum\ve_{T(i,j)}$, where we sum over the $(i,j)$ in the shape of $\lambda$. So the $i^{\rm th}$ component of the content of $T$ is the number of times that $i$ occurs in $T$. If $T$ is a tableau of shape $\lambda$ and weight $\mu$ with no repeated entries in the columns, then it is elementary to check that $\lambda^{\rm rev}\le\mu\le\lambda$, where $\lambda^{\rm rev}$ denotes the reversed tuple of $\lambda$.
We define the {\it weight} of a bitableau $(S,T)$ to be $(-\mu,\nu)\in X\times X$, where $\mu$ is the weight of $S$ and $\nu$ is the weight of $T$. If $H_{\mb Z}$ is the maximal torus of diagonal matrices of $\GL_{n,\mb Z}$, then the bideterminant $(S\,|\,T)$ is an $H_{\mb Z}\times H_{\mb Z}$ weight vector with weight equal to that of $(S,T)$.
The {\it degree} of a tableau or bitableau of shape $\lambda$ is defined to be $|\lambda|$. Note that the degree of a tableau is also determined by its content and that the degree of a bideterminant $(S\,|\,T)$ is equal to that of the bitableau $(S,T)$.

For a partition $\lambda$ of length $\le n$ we define $A_{\le\lambda}$ and $A_{<\lambda}$ to be the $\mb Z$-span of the bideterminants $(S\,|\,T)$ with $S$ and $T$ tableaux of shape $\le\lambda$ resp. $<\lambda$. Furthermore we define $\nabla_{\mb Z}(\lambda)$ and $\tilde\nabla_{\mb Z}(\lambda)$ as the $\mb Z$-span of the bideterminants $(T_\lambda\,|\,S)$ resp. $(S\,|\,T_\lambda)$ with $S$ a tableau of shape $\lambda$. Note that $A_{\le\lambda}$ and $A_{<\lambda}$ are $\GL_{n,\mb Z}\times\GL_{n,\mb Z}$-submodules of $\mb Z[\Mat_{n,\mb Z}]$ and that $\nabla_{\mb Z}(\lambda)$ and $\tilde\nabla_{\mb Z}(\lambda)$ are submodules for the right resp. left regular action. Note also that the notion of standardness is not involved in the above definitions. 

We now state a result from \cite{DeCEP}. {\it The reader should bear in mind that the transpose is (dominance) order reversing on partitions of the same number} (see e.g. \cite[1.11]{Mac}).
Assertion (i) below is due to Hesselink and Stein independently. It is an improved version of the first version of the straightening algorithm that can be found in \cite{DRS} and \cite{DKR}. The point is that one can show that the new shapes that show up during straightening are all $\le$ the original shape in the dominance order. As a consequence the combinatorial results on straightening match up nicely with the general theory of reductive groups. Of course, statements like (iv) below are known to hold for the coordinate rings of arbitrary reductive groups (see \cite[II.4.20]{Jan}), but the main point here is that these filtration subspaces and induced modules can be realised explicitly using bideterminants.

\begin{thm}[{\cite[Sect.~1-3]{DeCEP}}]\label{thm.GLn}
Let $\lambda$ be a partition of length $\le n$. Recall that the type $A_{n-1}$ partial order on weights (dominance order) is denoted by $\le$.
\begin{enumerate}[{\rm (i)}]
\item Let $S$ and $T$ be tableaux of shape $\lambda$. Then the bideterminant $(S\,|\,T)$ can be written as a linear combination $\sum_ia_i(S_i\,|\,T_i)$, where the $a_i$ are integers and the $S_i$ and $T_i$ are standard of shape $\le\lambda$ with the same content as $S$ and $T$ respectively.
\item The bideterminants $(S\,|\,T)$ with $S$ and $T$ standard form a basis of $\mb Z[\Mat_{n,\mb Z}]$.
\item The elements $(T_\lambda\,|\,T)$, $T$ standard of shape $\lambda$ form a basis of $\nabla_{\mb Z}(\lambda)$ and the elements $(T\,|\,T_\lambda)$, $T$ standard of shape $\lambda$ form a basis of $\tilde\nabla_{\mb Z}(\lambda)$.
\item The map $(S\,|\,T_\lambda)\ot(T_\lambda\,|\,T)\mapsto (S\,|\,T)$ defines an isomorphism
$$\tilde\nabla_{\mb Z}(\lambda)\ot_{\mb Z}\nabla_{\mb Z}(\lambda)\stackrel{\sim}{\to}A_{\le\lambda}/A_{<\lambda}$$
of $\GL_{n,\mb Z}\times\GL_{n,\mb Z}$-modules.
\end{enumerate}
\end{thm}

We recall some definitions from \cite{Don1} (see also \cite{Don2}). A subset $\pi$ of $X^+$ is called {\it saturated} if $\mu\in X^+$ and $\mu\le\lambda\in\pi$ implies $\mu\in\pi$. Now let $\pi$ be a saturated subset of $X^+$ and let $R$ be a principal ideal domain. For any torsion-free $\GL_{n,R}$ module (i.e. right $R[\GL_{n,R}]$-comodule) $M$ the submodule $O_\pi(M)$ is defined to be the sum (or the union) of the submodules of $M$ which are finitely generated (and therefore free of finite rank) over $R$ and whose dominant weights relative to the maximal torus of diagonal matrices lie in $\pi$. Clearly, $O_\pi(M)$ is the sum of the $O_{\pi'}(M)$, $\pi'\subseteq\pi$ finite, and $M/O_\pi(M)$ is torsion-free. Furthermore, $O_\pi(F\ot_{R}M)\cap M=O_\pi(M)$, where $F$ denotes the field of fractions of $R$. When we write $O_\pi(R[\Mat_{n,R}])$ or $O_\pi(R[\GL_{n,R}])$, then we consider $R[\Mat_{n,R}]$ or $R[\GL_{n,R}]$ as a $\GL_{n,R}$-module under the right regular action. The resulting module is stable under the $\GL_{n,R}\times\GL_{n,R}$-action.

Let $B_{\mb Z}$ and $B_{\mb Z}^-$ be the Borel subgroups of upper resp. lower triangular matrices in $\GL_{n,\mb Z}$, let $R$ be a commutative ring and let $\lambda\in X$. By applying $\lambda$ to the diagonal part of a upper or lower triangular matrix we obtain a $1$-dimensional representation of $B_{\mb Z}$ or $B_{\mb Z}^-$ which we also will denote by $\lambda$. Now let $H_\mb Z$ be $B_{\mb Z}$ or $B_{\mb Z}^-$. Then $\ind_{H_R}^{\GL_{n,R}}(\lambda)$ is defined to be the set of $f\in R[\GL_{n,R}]$ such that $f(hg)=\lambda(h)f(g)$ for all $h\in H(A)$ and $g\in\GL_n(A)$ and $A$ any commutative $R$-algebra. This is a submodule of $R[\GL_{n,R}]$ for the right regular action. In \cite{Jan} the induced module associated to $\lambda$ is defined by the property $f(gh)=\lambda(h)^{-1}f(g)$. This is a submodule for the left regular action. The automorphism of $R[\GL_{n,R}]$ given by the inversion maps one induced module onto the other and induces an isomorphism of $\GL_{n,R}$-modules. We have $\ind_{H_R}^{\GL_{n,R}}(\lambda)\cong R\ot_{\mb Z}\ind_{H_{\mb Z}}^{\GL_{n,\mb Z}}(\lambda)$, see \cite[II.8.8(1)]{Jan}. Furthermore, $\ind_{B_R^-}^{\GL_{n,R}}(\lambda)\cong \ind_{B_R}^{\GL_{n,R}}(\lambda^{\rm rev})$, where the isomorphism is given by left multiplication with the matrix of the permutation of $\un n$ that sends $i$ to $n+1-i$. Finally, we point out that the above definitions and facts about $O_\pi$ and the induced modules apply, with appropriate modifications, to any split reductive group scheme over $\mb Z$.

\begin{cornn}\
\begin{enumerate}[{\rm (i)}]
\item Let $\lambda\in\Lambda^+(n)$. Then $\nabla_{\mb Z}(\lambda)=\ind_{B_{\mb Z}^-}^{\GL_{n,\mb Z}}(\lambda)$ and $\tilde\nabla_{\mb Z}(\lambda)$ is the induced module of the $B_{\mb Z}$-module $-\lambda$, according to \cite{Jan}.
\item Let $\pi\subseteq\Lambda^+(n)$ be saturated. Then $O_\pi(\mb Z[\Mat_{n,\mb Z}])=O_\pi(\mb Z[\GL_{n,\mb Z}])$ is spanned by the bideterminants $(S\,|\,T)$ where $S$ and $T$ have shape $\in \pi$. Moreover, the $(S\,|\,T)$ with $S$ and $T$ standard of shape $\in \pi$ form a basis.
\end{enumerate}
\end{cornn}

\begin{proof}
(i).\ Over a field this is of course standard, see e.g. \cite[4.8]{Green}. We give a brief proof using the theory of reductive groups. One easily checks that $\nabla_{\mb Z}(\lambda)\subseteq\ind_{B_{\mb Z}^-}^{\GL_{n,\mb Z}}(\lambda)$ and that $\tilde\nabla_{\mb Z}(\lambda)$ is contained in the induced module of the $B_{\mb Z}$-module $-\lambda$, according to \cite{Jan}, which is isomorphic to $\ind_{B_{\mb Z}^-}^{\GL_{n,\mb Z}}(\lambda^*)$, where $\lambda^*=-\lambda^{\rm rev}$. Now $\tilde\nabla_{\mb Q}(\lambda)$ and $\nabla_{\mb Q}(\lambda)$ have dimension equal to the number of standard tableaux of shape $\lambda$ which is well known to be the dimension of the irreducible $\GL_n(\mb C)$-module of highest weight $\lambda$ and of course also that of its dual which has highest weight $\lambda^*$. The formal characters of these modules are given by Weyl's character formulas for $\lambda$ and $\lambda^*$. The formal character of $\ind_{B_{\mb Q}^-}^{\GL_{n,\mb Q}}(\mu)$ is given by Weyl's character formula for $\mu$, see \cite[Cor.II.5.11]{Jan}. So $\nabla_{\mb Q}(\lambda)=\ind_{B_{\mb Q}^-}^{\GL_{n,\mb Q}}(\lambda)$ and similarly for $\tilde\nabla_{\mb Q}(\lambda)$. Now the assertion follows from the fact that $\mb Z[\Mat_{n,\mb Z}]/\nabla_{\mb Z}(\lambda)$ and $\mb Z[\Mat_{n,\mb Z}]/\tilde\nabla_{\mb Z}(\lambda)$ are torsion-free (if $\preceq$ is the natural order, then $T_\lambda$ is standard).\\
(ii).\ We may assume that $\pi$ is finite. Let $A_\pi$ be the $\mb Z$-span of the bideterminants $(S\,|\,T)$ where $S$ and $T$ have shape $\in \pi$. Clearly $A_\pi\subseteq O_\pi(\mb Z[\Mat_{n,\mb Z}])\subseteq O_\pi(\mb Z[\GL_{n,\mb Z}])$. Furthermore, $\mb Q\ot_{\mb Z}A_\pi$ and $O_\pi(\mb Q[\GL_{n,\mb Q}])$ have the same dimension by \cite[3.2]{Don1} and the remarks above. So the assertion now follows from the fact that $\mb Z[\Mat_{n,\mb Z}]/A_\lambda$ and $O_\pi(\mb Z[\GL_{n,\mb Z}])/\mb Z[\Mat_{n,\mb Z}]$ are torsion-free.
\end{proof}

\begin{rems}\label{rems.GL}
1.~Theorem~\ref{thm.GLn} in \cite{DeCEP} was only proved for $\preceq$ the natural ordering. But one can in fact deduce from it a version for two different linear orderings $\preceq_1$ and $\preceq_2$ requiring in a bitableau $(S,T)$, $S$ to be standard relative to $\preceq_1$ and $T$ to be standard relative to $\preceq_2$. One simply has to apply the automorphism of the $\mb Z$-module $\mb Z[\Mat_{n,\mb Z}]$ induced by the automorphism $A\mapsto PAQ$ of the scheme $\Mat_{n,\mb Z}$ for suitable permutation matrices $P$ and $Q$. To obtain Theorem~\ref{thm.GLn}(i) and (ii) we take $P=Q^{-1}$ and to obtain assertion (iii) we first take $P=I$, the identity matrix, and then $Q=I$.\\
1.~If one works with nonsquare $n\times m$-matrices as in \cite{DeCEP}, then, in a bitableau $(S,T)$, $S$ should have entries in $\{1,\ldots,n\}$ and $T$ should have entries in $\{1,\ldots,m\}$. Furthermore, one has to work with $\Lambda^+(\min(n,m))$: partitions of length $\le\min(n,m)$.\\
2.~Let $\lambda_1,\lambda_2,\lambda_3,\ldots$ be an enumeration of the partitions of length $\le n$ such that $\lambda_i<\lambda_j$ implies $i<j$. Put $B_i=O_{\pi_i}(\mb Z[\Mat_{n,\mb Z}])$, where $\pi_i=\{\lambda_j\,|\,1\le j\le i\}$, a saturated set. Then $(B_i)_{i\ge0}$ is a $\GL_{n,\mb Z}\times\GL_{n,\mb Z}$-module filtration of $\mb Z[\Mat_{n,\mb Z}]$ with $B_i/B_{i-1}\cong\tilde\nabla_{\mb Z}(\lambda_i)\ot_{\mb Z}\nabla_{\mb Z}(\lambda_i)$ for $i\ge1$.
\end{rems}

We remind the reader of our convention to denote $\Mat_n(k)$ by $\Mat_n$. Clearly, the canonical epimorphism $k[(\Mat_n)_k]\to k[\Mat_n]$ is an isomorphism. The same remarks apply to $\GL_n$. In the next sections we will need some results relating graded pieces of the coordinate ring of a submonoid of $\Mat_n$ or $\Mat_l\times\Mat_m$ and the enveloping algebra of that monoid in a certain module. The result for submonoids of $\Mat_n$ is \cite[Prop.~1]{T} and its corollary. We now give the analogue for submonoids of $\Mat_l\times\Mat_m$.

Let $l,m$ be positive integers. Put $V=k^l$ and $W=k^m$. Any $(u,v)\in\Mat_l\times\Mat_m$ determines an endomorphism of $V^{\otimes r}\otimes W^{\otimes s}$ by
$$(u,v)(x_1\otimes\cdots\otimes x_r\otimes y_1\otimes\cdots\otimes y_s)=u(x_1)\otimes\cdots\otimes u(x_r)\otimes v(y_1)\otimes\cdots\otimes v(y_s).$$
For a subset $S$ of $\Mat_l\times\Mat_m$ we denote by $\mc{E}^{r,s}(S)$ the {\it enveloping algebra} of $S$ in $\End_k(V^{\otimes r}\otimes W^{\otimes s})$, that is, the subalgebra generated by the endomorphisms of $V^{\otimes r}\otimes W^{\otimes s}$ corresponding to the elements of $S$.
Using the isomorphism $\End_k(V^{\otimes r}\otimes W^{\otimes s})\cong\End_k(V^{\otimes r})\otimes\End_k(W^{\otimes s})$ we have
$\mc{E}^{r,s}(S)=\mc{E}^{r,0}(S)\otimes\mc{E}^{0,s}(S)$, where $\mc{E}^{r,0}(S)$ and $\mc{E}^{0,s}(S)$ are the enveloping algebras of $S$ in $\End_k(V^{\otimes r})$ and $\End_k(W^{\otimes s})$ respectively.

The algebra $k[\Mat_l\times\Mat_m]=k[\Mat_l]\otimes k[\Mat_m]$ is $\mb Z\times\mb Z$-graded. A subspace is homogeneous with respect to this grading if and only if it is stable under the action of $k^{\stimes}\times k^{\stimes}$ on $k[\Mat_l\times\Mat_m]$ which comes from the action on $\Mat_l\times\Mat_m$ given by $(a,b)\cdot (A,B)=(a A,b B)$. For any $\mb Z\times\mb Z$-graded vector space $U$ over $k$ we denote the graded piece of degree $(r,s)$ by $U^{r,s}$. For $S\subseteq\Mat_l\times\Mat_m$ we denote by $k[S]$ the $k$-algebra of $k$-valued functions on $S$ that are restrictions of functions in $k[\Mat_l\times\Mat_m]$. Clearly this notation conflicts with our earlier notation, e.g. in the case $S=\GL_l\times\GL_m$, so later on we will only use this notation in a situation where there is no conflict.

The next proposition and its corollary are a version of \cite[Prop.~1]{T} and its corollary. The proofs are a straightforward modification of the proofs in \cite{T}. For assertion (ii) has to use the fact that for groups $G_1$ and $G_2$, $U_1$ a $kG_1$-module and $U_2$ a $kG_2$-module, we have for the invariants: $(U_1\ot U_2)^{G_1\times G_2}=U_1^{G_1}\ot U_2^{G_2}$. The natural map from assertion (i) is given by precomposing with the natural homomorphism $M\to\mc{E}^{r,s}(M)$ of monoids.


\begin{prop}\label{prop.envalg}
Let $M$ be a submonoid of $\Mat_l\times\Mat_m$ with $(k^{\stimes}\times k^{\stimes})M=M$. Then
\begin{enumerate}[{\rm(i)}]
\item the natural map $\mc{E}^{r,s}(M)^*\to k[M]^{r,s}$ is an isomorphism of coalgebras,
\item $\mc{E}^{r,s}(\GL_l\times\GL_m)=\End_{k\la\Sym_r\times\Sym_s\ra}(V^{\otimes r}\otimes W^{\otimes s})$.
\end{enumerate}
\end{prop}

\begin{cornn}
Let $M$ be a submonoid of $\Mat_l\times\Mat_m$ with $(k^{\stimes}\times k^{\stimes})M=M$, let $I$ be the (homogeneous) ideal of polynomial functions on $\Mat_l\times\Mat_m$ that vanish on $M$. Furthermore, let $g_1,\ldots,g_t$ be nonzero homogeneous elements of $I$. Denote the isomorphism $k[\Mat_l\times\Mat_m]^{r,s}\to \mc{E}^{r,s}(\Mat_l\times\Mat_m)^*$ by $\eta$. Then
the elements $g_1,\ldots,g_t$ are generators of $I$ if and only if for each $r,s\ge 0$, $(r,s)\ne(0,0)$, the functionals $\eta(g_im_i)$, where the $m_i$ are arbitrary monomials in the matrix entries of degree $(r,s)-\deg(g_i)$, define the algebra $\mc{E}^{r,s}(M)$.
\end{cornn}

\section{Symplectic straightening}\label{s.symplecticstraightening}
From now on we assume that $n=2m$, $m\ge1$, is even.
We take
\begin{equation}\label{eq.J}
J=\begin{bmatrix}
0&J_m\\
-J_m&0
\end{bmatrix}
\text{ or }
J=\begin{bmatrix}
0&I\\
-I&0
\end{bmatrix},
\end{equation}
where $J_m$ is the $m\times m$ matrix with ones on the antidiagonal and zeros elsewhere, and $I$ is the $m\times m$ identity matrix. Everything in this section will be valid for both choices of $J$. For every integer $i\in\un n$ there is a unique integer $i'\in\un n$ and a unique nonzero scalar $\e_i\in k$ such that $J_{i,i'}=\e_i$. Clearly $i''=i$ and we have $\e_i=1$ if $i\le m$ and $\e_i=-1$ if $i>m$.

Let $V_{\mb Z}=\mb Z^n$ be the natural module of the monoid scheme $\Mat_{n,\mb Z}$ and denote the standard basis elements by $v_1,\dots,v_n$. On $V_{\mb Z}$ we define the nondegenerate symplectic form $\la\ ,\ \ra$ by
$$\la x,y\ra:=x^TJy=\sum_{i=1}^n\e_ix_iy_{i'}\ .$$
Let $R$ be a commutative ring. The {\it symplectic group} $\Sp_n(R)$ over $R$ consists of the $n\times n$-matrices $A$ over $R$ that satisfy $A^TJA=J$, i.e. the matrices for which the corresponding automorphism of $V$ preserves the form $\la\ ,\ \ra$.
The {\it symplectic monoid} $\SpM_n(R)$ over $R$ is defined as the set of $n\times n$-matrices $A$ over $R$ for which there exists a scalar $d(A)\in R$ such that $A^TJA=AJA^T=d(A)J$. The group of invertible elements of $\SpM_n(R)$ is the {\it symplectic similitude group} $\GSp_n(R)$ over $R$. It consists of the matrices $A$ that satisfy $A^TJA=d(A)J$ for some invertible scalar $d(A)\in R$, i.e. the invertible matrices for which the corresponding automorphism of $V$ preserves the form $\la\ ,\ \ra$ up to a scalar.
We denote the functors $R\mapsto\Sp_n(R)$, $R\mapsto\SpM_n(R)$ and $R\mapsto\GSp_n(R)$ by $\Sp_{n,\mb Z}$, $\SpM_{n,\mb Z}$ and $\GSp_{n,\mb Z}$. The functors $\Sp_{n,\mb Z}$ and $\SpM_{n,\mb Z}$ are closed subschemes of $\Mat_{n,\mb Z}$ and the functor $\GSp_{n,\mb Z}$ is a closed subscheme of $\GL_{n,\mb Z}$. For $i,j\in\un n$, define $g_{ij},\ov{g}_{ij}\in\mb Z[\Mat_{n,\mb Z}]$ by
\begin{equation}\label{eq.g}
g_{ij}:=\sum_{l=1}^n\e_lx_{li}x_{l'j}\text{\quad and\quad}\ov{g}_{ij}:=\sum_{l=1}^n\e_lx_{il}x_{jl'}\ .
\end{equation}
Note that $g_{ii}=\ov g_{ii}=0$ and that $g_{ij}=-g_{ji}$ and $\ov g_{ij}=-\ov g_{ji}$.
The ideal of $\Sp_{n,\mb Z}$ in $\mb Z[\Mat_{n,\mb Z}]$ is generated by the elements $g_{ij}$, $1\le i<j\le n$, $i\ne j'$, and $g_{rr'}-1$, $1\le r\le m$.
The ideal of $\GSp_{n,\mb Z}$ in $\mb Z[\GL_{n,\mb Z}]$ is generated by the elements $g_{ij}$, $1\le i<j\le n$, $i\ne j'$, and $g_{rr'}-g_{ss'}$, $1\le r,s\le m$.
The ideal of $\SpM_{n,\mb Z}$ in $\mb Z[\Mat_{n,\mb Z}]$ is generated by the elements
\begin{equation}\label{eq.SpM}
\{g_{ij},\,\ov g_{ij},\,g_{rr'}-\ov g_{ss'}\,|\,1\le i<j\le n, i\ne j', 1\le r,s\le m\}.
\end{equation}
Using \cite[Thm.~II.5.2.1]{DemGa} or \cite[12.2]{W}, a simple Lie algebra computation shows that the fibers of $\Sp_{n,\mb Z}$ and $\GSp_{n,\mb Z}$ over $\Spec(\mb Z)$ are reduced. This means that for $F$ a prime field (including $\mb Q$) the algebras $F[\Sp_{n,F}]$ and $F[\GSp_{n,F}]$ are reduced. Since these fields are perfect this holds for any field (as one could have showed directly by the same method). So $K[\Sp_{n,K}]$ and $K[\GSp_{n,K}]$ are the coordinate rings of the connected reductive algebraic groups $\Sp_n(K)$ and $\GSp_n(K)$, see \cite{Bo}. Furthermore it is clear that they are defined over the prime field as closed subgroups of $\GL_n(K)$. The derived group of $\GSp_n(K)$ is $\Sp_n(K)$. Clearly $\GSp_n=\GSp_n(k)$ is the group of $k$-points of $\GSp_n(K)$ and therefore it is dense in $\GSp_n(K)$ by \cite[Cor.~V.18.3]{Bo}.
Put $d=g_{11'}$. The restriction of $d$ to $\SpM_{n,\mb Z}$ is called the {\it coefficient of dilation}. It is equal to the function $d$ mentioned above. We have $d^n=\det^2$ in $\mb Z[\SpM_{n,\mb Z}]$ (after Theorem~\ref{thm.SpM} one can show that $d^m=\det$). 
Note that $\mb Z[\GSp_{n,\mb Z}]=\mb Z[\SpM_{n,\mb Z}][d^{-1}]$. Of course there is, just as in the $\GL_{n,\mb Z}$-case, an action of $\Sp_{n,\mb Z}\times\Sp_{n,\mb Z}$ on $\mb Z[\SpM_{n,\mb Z}]$, $\mb Z[\GSp_{n,\mb Z}]$ and $\mb Z[\Sp_{n,\mb Z}]$.

We denote the character group of the maximal torus of $\SpM_{n,\mb Z}$ of diagonal matrices by $X$. We have $X\cong\mb Z^m$ and we denote the standard basis elements by $\ve_i$. We embed $Z^m$ in $\mb Z^n$ by extending $m$-tuples with $m$ zeros. The restriction of a character of the diagonal matrices in $\GL_{n,\mb Z}$ to those in $\Sp_{n,\mb Z}$ is given by the map $\lambda\mapsto\ov\lambda:Z^n\to\mb Z^m$ with $\ov\ve_i=\ve_i$ for $i\le m$ and $\ov\ve_i=-\ve_i$ for $i>m$. In the root system of type $C_m$ we choose the set of positive roots as usual, they are: $\ve_i-\ve_j$, $1\le i<j\le m$, $\ve_i+\ve_j$, $i,j\in\un m$ with $i\ne j$ and $2\ve_i$, $i\in\un m$. If $J$ equals the first matrix in \eqref{eq.J}, then the corresponding Borel subgroup of $\Sp_{n,\mb Z}$ is the subgroup of upper triangular matrices. If $J$ equals the second matrix in \eqref{eq.J}, then the corresponding Borel subgroup of $\Sp_{n,\mb Z}$ is the subgroup of matrices of the form $
\begin{bmatrix}
A&B\\
0&C
\end{bmatrix}
$, with $A$ upper triangular and $C$ lower triangular. The set of dominant weights is $X^+=\Lambda^+(m)$, the partitions of length $\le m$. We denote the type $C_m$ partial order on $\mb Z^m$ by $\le$. We have $\mu\le\lambda$ if and only if $|\lambda|-|\mu|$ is even $\ge0$ and $\sum_{i=1}^j\mu_i\le\sum_{i=1}^j\lambda_i$ for all $j\in\un m$.

We assume given a linear order $\preceq$ on $\un n$ such that for all $i\in \un m$, $i$ is the immediate successor of $i'$ or the other way around, for example $1'\prec1\prec2'\prec2\cdots\prec m'\prec m$. Define $\zeta:\un n\to\un m$ by $|\{i\in\un n\,|\,i\preceq j\text{ or }i\preceq j'\}|=2\zeta(j)$ for all $j\in\un n$. Note that $i\preceq j$ implies $\zeta(i)\le\zeta(j)$ and that $\zeta(i)=\zeta(j)$ if and only if $i=j$ or $j'$. A subset $I$ of $\un n$ is called {\it symplectic standard} if $$|\{i\in I\,|\,i\preceq j\text{ or }i\preceq j'\}|\le\zeta(j)$$ for all $j\in\un m$. Taking $j\in \un m$ such that either $j$ or $j'$ is the maximal element of $I$ according to $\preceq$, we see that $|I|\le m$ whenever $I$ is symplectic standard. We identify each subset $I$ of $\un n$ with the one column tableau whose entries are the elements of $I$ and whose entries are strictly increasing (according to $\preceq$) from top to bottom. A tableau is called {\it symplectic standard} if it is $\GL_n$-standard (relative to $\preceq$) and if the first column is symplectic standard as a set. So the shape of a symplectic standard tableau has length $\le m$. If a tableau is symplectic standard, then all its columns are symplectic standard. So a $\GL_n$-standard tableau is symplectic standard if and only if for all $i\in\un m$ the occurrences of $i$ and $i'$ are limited to the first $\zeta(i)$ rows.
Note that for two linear orderings $\preceq_1$ and $\preceq_2$ of $\un n$ as above there is a permutation of $\un n$ which stabilises $\{\{i,i'\}\,|\,i\in\un m\}$ and induces a bijection between the tableaux that are symplectic standard relative to $\preceq_1$ and those that are symplectic standard relative to $\preceq_2$.
A bitableau $(S,T)$ is called symplectic standard if $S$ and $T$ are symplectic standard. We denote the restriction of a bideterminant $(S\,|\,T)$ to $\SpM_{n,\mb Z}$ by the same symbol. The {\it symplectic weight} of a tableau with $\GL_n$ weight $\mu$ is defined as the restriction $\ov\mu$. If $T$ is a tableau with shape $\lambda$ of length $\le m$, with symplectic weight $\mu$ and with no repeated entries in the columns, then it is easy to check that $-\lambda\le\mu\le\lambda$. The symplectic weight of a bitableau is also defined by restriction of its $\GL_n$ weight.

Let $\bigwedge V_{\mb Z}$ be the exterior algebra on $V_{\mb Z}$. We denote the set of $r$ element subsets of $\un n$ by $P(n,r)$. For $I=\{i_1,\ldots,i_r\}\subseteq\un n$ with $i_1\prec i_2\prec\cdots\prec i_r$ we define $v_I:=v_{i_1}\wedge\cdots\wedge v_{i_r}$ and if $I\subseteq\un m$, then we define $z(I)=v_{i_1'}\wedge v_{i_1}\wedge\cdots\wedge v_{i_r'}\wedge v_{i_r}$. Recall that an element of the exterior algebra whose odd degree components are $0$ is central. For $r\in\{0,\ldots,m\}$, we define $z_r=\sum_{I\in P(m,r)}z(I)$ ($z_0=1$). Note that $z_1\in\bigwedge^2 V_{\mb Z}$ is the element corresponding to the symplectic form under the canonical isomorphism $\bigwedge^2 V_{\mb Z}\cong\bigwedge^2 V_{\mb Z}^*$.
As is well-known, the $v_I$ form a basis of $\bigwedge V_{\mb Z}$. In $\bigwedge V_{\mb Q}$ we have $z_r=(1/r!)z_1^r$. From this we deduce $z_rz_s=\binom{r+s}{r}z_{r+s}$.

The {\it symplectic content} of a tableau $T$ is the tuple $a\in\mb Z^m$, such that $a_i$ is the number of occurrences of $i$ and $i'$ in $T$. We define the lexicographical order $\trianglelefteq$ on $\mb Z^m$ as follows: $a\trianglelefteq b$ if $a=b$ or $a\ne b$ and $a_i<b_i$, where $i=\max\{j\in\un m\,|\,a_j\ne b_j\}$ and the maximum is taken according to $\preceq$. For subsets $I$ and $J$ with symplectic contents $a$ and $b$ respectively, we write $I\trianglelefteq J$ if $a\trianglelefteq b$.

\begin{prop}[cf. \cite{Ber},\cite{Don4},\cite{Oe}]\label{prop.extbasis}\
\begin{enumerate}[{\rm (i)}]
\item For every $J\in P(n,r)$ not symplectic standard, the vector $v_J\in\bigwedge^r V_{\mb Z}$ can be written as
\begin{equation}
v_J=\sum_L a_{JL}v_L+\sum_{t,L} b_{JL}z_t\wedge v_L,\label{eq.v}
\end{equation}
with $a_{JL},b_{JL}\in\mb Z$; the first sum over all $L\in P(n,r)$ symplectic standard
with $L\vartriangleright J$ and the second sum over all $t\in\{1,\ldots,\lfloor r/2\rfloor\}$ and $L\in P(n,r-2t)$ symplectic standard. Furthermore, all the $L$ occurring have the same symplectic weight as $J$.
\item The vectors $z_t\wedge v_J$, $t\in\{0,\ldots,\lfloor r/2\rfloor\}$, $J\in P(n,r-2t)$ symplectic standard, form a basis of $\bigwedge^r V_{\mb Z}$.
\end{enumerate}
\end{prop}

\begin{proof}
(i).\ Let $J\in P(n,r)$ be not symplectic standard. To prove the first assertion it suffices to show that $v_J$ can be written as in \eqref{eq.v} with the first sum over all $L\in P(n,r)$ with $L\vartriangleright J$ and the second sum over all $t\in\{1,\ldots,\lfloor r/2\rfloor\}$ and $L\in P(n,r-2t)$, since then we can finish by induction on $|J|$ and $\trianglelefteq$. Now any element of the $r^{\rm th}$ graded piece of the ideal of $\bigwedge^r V_{\mb Z}$ generated by $z_1,\ldots,z_m$ can be written as in the second sum. So what we need is precisely what is proved in \cite[Lemma~8.1]{Oe}.
Let $H_{\mb Z}$ be the maximal torus of diagonal matrices in $\Sp_{n,\mb Z}$. Then the final assertion follows by applying the projection onto the $H_{\mb Z}\times H_{\mb Z}$ weight space to which $v_J$ belongs to \eqref{eq.v}.\\
(ii).\ By \cite[Thm~17.5]{FH} $\bigwedge^r V_{\mb C}$ is the direct sum of the irreducible $\Sp_n(\mb C)$ representations of highest weights $r\ve_1, (r-2)\ve_1,\ldots$. By \cite{King} the dimension of the irreducible representation of highest weight $(r-2t)\ve_1$ is equal to the number of symplectic standard $J\in P(n,r-2t)$. So the assertion follows from (i) and the fact that the canonical map $\bigwedge^r V_{\mb Z}\to\bigwedge^r V_{\mb C}$ is an embedding, since $\bigwedge^r V_{\mb Z}$ is free.
\end{proof}

For the comodule map $\Delta_\wedge:\bigwedge^r V_{\mb Z}\to\bigwedge^r V_{\mb Z}\ot\mb Z[\SpM_{n,\mb Z}]$ of the $\SpM_{n,\mb Z}$-action we have
\begin{equation}\label{eq.comod1}
\Delta_\wedge(v_J)=\sum_Iv_I\ot(I\,|\,J)\, ,
\end{equation}
where the sum is over all $I\in P(n,r)$. This just follows from the corresponding equations for the comodule map of the $\Mat_{n,\mb Z}$-action by restriction. Note that $\Delta_\wedge$ is a homomorphism of algebras, since $\Mat_{n,\mb Z}$ acts on $\bigwedge V_{\mb Z}$ by algebra endomorphisms. We record now a result from \cite{Oe}. For $t=1$ the result follows from the relations \eqref{eq.SpM}. As pointed out in \cite{Oe} the result would follow immediately from the equality $z_t=(1/t!)z_1^t$ in $\bigwedge V_{\mb Q}$ if we would know that $\mb Z[\SpM_{n,\mb Z}]$ has no torsion. This follows immediately from Theorem~\ref{thm.SpM}, but to prove that theorem we need the result below.

\begin{lem}[\cite{Oe}]\label{lem.semiinv}
For each $t\in\un m$, the element $z_t\in\bigwedge^{2t} V_{\mb Z}$ is a semi-invariant of $\SpM_{n,\mb Z}$ with weight $d^t$, that is,
$\Delta_\wedge(z_t)=z_t\ot d^t$.
\end{lem}

Define the integers $c_{LI}$, $I\in P(n,r)$, $L\in P(n,r-2t)$,  by $z_t\wedge v_L=\sum_Ic_{LI}v_I$. Then we get the following corollary to Proposition~\ref{prop.extbasis}.

\begin{corgl}
Let $I,J\in P(n,r)$ with $J$ not symplectic standard. Then we have in $\mb Z[\SpM_{n,\mb Z}]$
\begin{align}
(I\,|\,J)=\sum_L a_{JL}(I\,|\,L)+\sum_{t,L,L'} b_{J,L}c_{L'I}d^t(L'\,|\,L)&\text{\quad and}\label{eq.tab1}\\
(J\,|\,I)=\sum_L a_{JL}(L\,|\,I)+\sum_{t,L,L'} b_{J,L}c_{L'I}d^t(L\,|\,L')&\label{eq.tab2}\, ,
\end{align}
with $a_{JL},b_{JL}\in\mb Z$ given by \eqref{eq.v}; in both cases the first is sum over all $L\in P(n,r)$ symplectic standard
with $L\vartriangleright J$ and the second sum is over all $t\in\{1,\ldots,\lfloor r/2\rfloor\}$ and $L,L'\in P(n,r-2t)$ symplectic standard. Furthermore, all the $L$ and $L'$ occurring have the same symplectic weight as $J$ and $I$ respectively.
\end{corgl}

\begin{proof}
Equation \eqref{eq.tab1} follows by applying the comodule map \eqref{eq.comod1} to \eqref{eq.v}, using the fact that $\Delta_\wedge$ is a homomorphism of algebras and using Proposition~\ref{lem.semiinv} and the definition of the $c_{LI}$. The automorphism $x_{ij}\mapsto x_{ji}$ of $\mb Z[\Mat_{n,\mb Z}]$ sends the bideterminant $(S\,|\,T)$ to $(T\,|\,S)$. Furthermore, it leaves the ideal generated by the elements \eqref{eq.g} stable, so it induces an automorphism $\varphi$ of $\mb Z[\SpM_{n,\mb Z}]$. Note that it is important here that both the relations corresponding to the condition $A^TJA=d(A)J$ and those corresponding to the condition $AJA^T=d(A)J$ occur in \eqref{eq.g}. Equation \eqref{eq.tab2} now follows by applying $\varphi$ to \eqref{eq.tab1}. Note that we could have obtained \eqref{eq.tab2} directly by considering the right action of $\SpM_{n,\mb Z}$ on $\bigwedge V_{\mb Z}$. Let $H_{\mb Z}$ be the maximal torus of diagonal matrices in $\Sp_{n,\mb Z}$. Then the final assertion follows by applying the projections onto the $H_{\mb Z}\times H_{\mb Z}$ weight spaces to which $(I\,|\,J)$ and $(J\,|\,I)$ belong to \eqref{eq.tab1} and \eqref{eq.tab2} respectively.
\end{proof}

\begin{rems}\label{rems.sympstraight}\ \\
1.\ If $r>m$, then no $L\in P(n,r)$ can be symplectic standard. So for $r>m$ and $J\in P(n,r)$ the first sums in \eqref{eq.v}, \eqref{eq.tab1} and \eqref{eq.tab2} are zero.\\
2.\ Let $L\subseteq\un n$ and $t\in\{1,\ldots,\lfloor r/2\rfloor\}$. Then $z_t\wedge v_L$ is a signed sum of $v_I$ with $I\cap I'\ne\emptyset$ ($I':=\{i'\,|\,i\in I\}$), so $c_{LI}=0$ for all $I\subseteq\un n$ with $I\cap I'=\emptyset$. So for such $I$ the second sums in \eqref{eq.tab1} and \eqref{eq.tab2} are zero. Note that this applies when $I\subseteq\un m$.\\
3.\ Let $N$ be the ideal of $\in\bigwedge V_{\mb Z}$ generated by $z_1,\ldots,z_m$. As pointed out in the proof of the above corollary one needs to prove straightening in $\in\bigwedge V_{\mb Z}/N$, since then one can prove it in general by induction on the degree. The relations that are used for this are $$z(I)\equiv(-1)^r\sum_{J\in P(m,r),J\cap I=\emptyset}z(J)\qquad({\rm mod}\ N)$$ for all $I\in P(m,r)$. See \cite[(19)]{Oe} or \cite[2.2]{Don4}. If we now multiply both sides with $v_L$, for each $L\subseteq\un n$ such that $I\cap(L\cup L')=\emptyset$, then we may restrict the sum to all $J\in P(m,r)$ such that $J\cap(I\cup L\cup L')=\emptyset$ and applying the comodule map we obtain explicit relations that can be used for the symplectic straightening. See condition (iv) on page 119 in \cite{Don4} and also \cite[Prop.~1.8]{DeC}.
\end{rems}

For the proof of Theorem~\ref{thm.SpM} we need the following simple lemma.

\begin{lem}\label{lem.dim}
Let $N$ be a finitely generated $\mb Z$-module and let $x_1,\ldots,x_r$ be generators of $N$.
If, for some field $\mb F$ of characteristic $0$, $\dim\mb F\otimes_{\mb Z}N\ge r$, then
$(x_1,\ldots,x_r)$ is a basis of $N$ (and, of course, $\dim\mb F\otimes_{\mb Z}N=r$).
\end{lem}

\begin{proof}
Assume that $N$ has torsion. Then we have for some prime $p$ that $\dim \mb F_p\otimes_{\mb Z}N>{\rm rk}\,N=\dim \mb F\otimes_{\mb Z}N\ge r$, where $\mb F_p$ denotes the field with $p$ elements. This contradicts the fact that $(1\ot x_1,\ldots,1\ot x_r)$ generates $\mb F_p\otimes_{\mb Z}N$. So the canonical map $N\to\mb F\otimes_{\mb Z}N$ is an embedding. Since $(x_1,\ldots,x_r)$ generates $\mb F\otimes_{\mb Z}N$, it must be a basis of $\mb F\otimes_{\mb Z}N$ and therefore also be independent in $N$.
\end{proof}
\vfill\eject
\begin{thm}[cf.\,{\cite[Thm.\,6.1]{Oe}}]\label{thm.SpM}
Recall that the type $C_m$ partial order on weights is denoted by $\le$.
\begin{enumerate}[{\rm (i)}]
\item Let $\lambda$ be a partition of length $\le m$ and let $S$ and $T$ be tableaux of shape $\lambda$. Then the bideterminant $(S\,|\,T)\in \mb Z[\SpM_{n,\mb Z}]$ can be written as a linear combination $\sum_ia_id^{t_i}(S_i\,|\,T_i)$, where the $a_i$ are integers and $(S_i,T_i)$ is symplectic standard of shape $\mu_i\le\lambda$ with $|\mu_i|=|\lambda|-2t_i$ and with the same symplectic weight as $(S,T)$.
\item The elements $d^t(S\,|\,T)$ with $t$ an integer $\ge0$ and $S$ and $T$ symplectic standard form a basis of $\mb Z[\SpM_{n,\mb Z}]$.
\end{enumerate}
\end{thm}
\begin{proof}
(i).\ Note that for $\lambda,\mu\in\mb Z^m$, we have that $\lambda\le_{C_n}\mu$ implies $\lambda\le\mu$ ($\le$ is $\le_{C_m}$). So it suffices to prove the result for $\lambda$ a partition of length $\le n$ and with the ordering $\le$ on weights replaced by $\le_{C_n}$. The point is that we will use $\GL_n$-straightening and there shapes of length $>m$ may show up (of course a tableau with such a shape can never be symplectic standard). First we show that, when $S$ and $T$ are $\GL_n$-standard, but $S$ or $T$ is not symplectic standard, $(S\,|\,T)$ can be written as a $\mb Z$-linear combination $\sum_ia_i(S_i\,|\,T_i)+\sum_ib_id^{t_i}(S'_i\,|\,T'_i)$, where in the first sum $S_i$ and $T_i$ are of shape $\lambda$ and $S_i$ or $T_i$ has symplectic content $\vartriangleright$ than that of $S$ or $T$ respectively and in the second sum the $t_i$ are $\ge1$ and the $S'_i$ and $T'_i$ are of shape $\mu_i\le_{C_n}\lambda$ with $|\mu_i|=|\lambda|-2t_i$. To prove this we apply \eqref{eq.tab1} or \eqref{eq.tab2} to the minor corresponding to the first columns of $S$ and $T$, and all that remains to check is that the new shapes that arise from the second sums in \eqref{eq.tab1} or \eqref{eq.tab2} are $\le_{C_n}\lambda$. Such a shape $\mu$ is formed by shortening the length of the first column of $\lambda$ by an even number, $2t$ say, and moving the resulting column to the right position to get a Young diagram. So $\mu=\lambda-\sum_{i=0}^{t-1}(\e_{l-2i}+\e_{l-1-2i})\le_{C_n}\lambda$, where $l$ is the length of $\lambda$.

Now one can finish by induction. The argument via contradiction is as follows. Give $\mb Z^m\times\mb Z^m$ the lexicographical ordering based on the ordering $\trianglelefteq$ of $\mb Z^m$. Assume the assertion doesn't hold. Pick a counterexample with $|\lambda|$ minimal and then with the pair of symplectic contents maximal with respect to the ordering we just defined. By Theorem~\ref{thm.GLn}(i) and the fact that $\mu\le_{A_{n-1}}\nu$ implies $\mu\le_{C_n}\nu$, we may assume that our bitableau is $\GL_n$-standard. Now what we proved above leads to a contradiction.\\
(ii).\ Clearly, the piece of degree $r$ of $\mb C[\SpM_{n,\mb C}]$ surjects onto that of $\mb C[\ov{\GSp_{n}(\mb C)}]$, where the bar denotes Zariski closure. By \cite[Prop.~1(i)]{T} this is the dual space of the enveloping algebra $E$ of $\Sp_{n}(\mb C)$ (or $\GSp_n(\mb C)$) in $\End_\mb C(V_\mb C^{\ot r})$. By Weyl's complete reducibility theorem for complex semisimple Lie algebras $E$ has dimension $\sum_\lambda \dim L(\lambda)^2$, where $L(\lambda)$ is the irreducible $\Sp_n(\mb C)$-module of highest weight $\lambda$ and the sum is over all $\lambda$ such that $L(\lambda)$ appears in $V_\mb C^{\ot r}$. By \cite[VI.3]{Wey}(p.175) these are the partitions of $r, r-2,\ldots$. Now $\dim L(\lambda)$ is the number of symplectic standard tableaux of shape $\lambda$ by \cite{King}, so $\dim E$ is the number of bitableaux $(S,T)$ with $S$ and $T$ symplectic standard of degree $r, r-2,\ldots$. So, by (i) and Lemma~\ref{lem.dim}, the elements $d^t(S\,|\,T)$ with $t\in\{0,\ldots,\lfloor r/2\rfloor\}$ and $S$ and $T$ symplectic standard of degree $r-2t$ form a basis of the degree $r$ piece of $\mb Z[\SpM_{n,\mb Z}]$.
\end{proof}\vfill

By Theorem~\ref{thm.SpM}(ii) over $K$, $d\in K[\SpM_{n,K}]$ is not a zero divisor. Therefore $K[\SpM_{n,K}]\subseteq K[\SpM_{n,K}][d^{-1}]=K[\GSp_n(K)]$ is a domain. So $\SpM_{n,K}$ is reduced and $\SpM_n(K)$ is an irreducible closed subvariety of $\Mat_n(K)$ which is defined over the prime field. Since $\GSp_n$ is dense in $\GSp_n(K)$, it is also dense in $\SpM_n(K)$, so the canonical epimorphism $k[(\SpM_n)_k]\to k[\SpM_n]$ is an isomorphism. This was deduced in a different way in \cite[Cor.~6.2]{Oe}. Note that the reducedness results above can be expressed by saying that fibers of $\SpM_{n,\mb Z}$ over $\Spec(\mb Z)$ are reduced. In case $k=K$ is algebraically closed, we have that $k[\SpM_n]$ is a unique factorisation domain by \cite[Prop.~2]{T}. In particular, it is normal and therefore Cohen-Macaulay by \cite[Cor~2]{Rit}.

For $\lambda$ a partition of length $\le m$ we define $A_{\le\lambda},A_{<\lambda},\nabla(\lambda),\tilde\nabla(\lambda)\subseteq\mb Z[\Sp_{n,\mb Z}]$ completely analogous to the $\Mat_{n,\mb Z}$-case. Note that the automorphism of $\mb Z[\Sp_{n,\mb Z}]$ given by $A\mapsto JA^TJ=A^{-1}$ on $\Sp_{n,\mb Z}$ maps $\nabla(\lambda)$ and $\tilde\nabla(\lambda)$ onto each other and defines an isomorphism between the two.
\begin{thm}\label{thm.Sp}
Let $\lambda$ be a partition of length $\le m$.
\begin{enumerate}[{\rm (i)}]
\item Let $S$ and $T$ be tableaux of shape $\lambda$. Then the bideterminant $(S\,|\,T)\in \mb Z[\Sp_{n,\mb Z}]$ can be written as a linear combination $\sum_ia_i(S_i\,|\,T_i)$, where the $a_i$ are integers and the $(S_i,T_i)$ are symplectic standard of shape $\le\lambda$ with the same symplectic weight as $(S,T)$.
\item The bideterminants $(S\,|\,T)$ with $S$ and $T$ symplectic standard form a basis of $\mb Z[\Sp_{n,\mb Z}]$.
\item The elements $(T_\lambda\,|\,T)$, $T$ symplectic standard of shape $\lambda$ form a basis of $\nabla_{\mb Z}(\lambda)$ and the elements $(T\,|\,T_\lambda)$, $T$ symplectic standard of shape $\lambda$ form a basis of $\tilde\nabla_{\mb Z}(\lambda)$.
\item The map $(S\,|\,T_\lambda)\otimes(T_\lambda\,|\,T)\mapsto (S\,|\,T)$ defines an isomorphism
$$\tilde\nabla_{\mb Z}(\lambda)\otimes_{\mb Z}\nabla_{\mb Z}(\lambda)\stackrel{\sim}{\to}A_{\le\lambda}/A_{<\lambda}$$
of $\Sp_{n,\mb Z}\times\Sp_{n,\mb Z}$-modules.
\end{enumerate}
\end{thm}

\begin{proof}
We have $\mb Z[\Sp_{n,\mb Z}]=\mb Z[\SpM_{n,\mb Z}]/(d-1)$, so (i) and (ii) follow from Theorem~\ref{thm.SpM}.\\
(iii).\ We first show that the given tableaux span $\nabla_{\mb Z}(\lambda)$ and $\tilde\nabla_{\mb Z}(\lambda)$. By Theorem~\ref{thm.GLn}(iii) we only have to show that a symplectic straightening step \eqref{eq.tab1} or \eqref{eq.tab2} applied to $(T_\lambda\,|\,T)$ or $(T\,|\,T_\lambda)$, $T$ not symplectic standard, of shape $\lambda$, does not yield any new shapes. This follows from Remark~\ref{rems.sympstraight}.2. If we take $\preceq$ to be the natural order, then $T_\lambda$ is standard, so (ii) gives us then that $\nabla_{\mb Z}(\lambda)$ and $\tilde\nabla_{\mb Z}(\lambda)$ are free of rank equal to the number of symplectic standard tableaux of shape $\lambda$. Since this number is independent of the linear order $\preceq$, this proves (iii).\\
(iv). The argument is precisely the same as in \cite{DeCEP}. It suffices to show that the map is well-defined, since then it will map a basis to a basis. This follows from the fact that, once we factor out all bideterminants of shapes $<\lambda$, straightening steps for $(S\,|\,T)$ arising from $S$ are valid for all $T$ and vice versa. For the symplectic straightening steps this is expressed by the fact that the integers $a_{JL}$ in \eqref{eq.tab1} and \eqref{eq.tab2} do not depend on $I$.
\end{proof}

Let $B_{\mb Z}$ and $B_{\mb Z}^-$ be the Borel subgroups in $\Sp_{n,\mb Z}$ corresponding to the positive resp. negative roots. The notion of saturated set, the functor $O_\pi$ and the induced modules are defined completely analogous to the $\GL_{n,\mb Z}$-case in Section~1. The proof of the corollary below is also completely analogous to the Corollary to Theorem~\ref{thm.GLn}, so we omit it. Assertion (i) was proved in \cite{Don4} over $K$ in another way.

\begin{cornn}\
\begin{enumerate}[{\rm (i)}]
\item Let $\lambda\in\Lambda^+(m)$. Then $\nabla_{\mb Z}(\lambda)$ is the induced module $\ind_{B_{\mb Z}^-}^{\Sp_{n,\mb Z}}(\lambda)$ and $\tilde\nabla_{\mb Z}(\lambda)$ is the induced module of the $B_{\mb Z}$-module $-\lambda$, according to \cite{Jan}.
\item Let $\pi\subseteq\Lambda^+(m)$ be saturated. Then $O_\pi(\mb Z[\Sp_{n,\mb Z}])$ is spanned by the bideterminants $(S\,|\,T)$ where $S$ and $T$ have shape $\in \pi$. Moreover, the $(S\,|\,T)$ with $S$ and $T$ symplectic standard of shape $\in \pi$ form a basis.
\end{enumerate}
\end{cornn}

\begin{rem}
The main results in \cite{DeC} are also valid for the symplectic standard tableaux of King that we used. First we note that if we take for $\preceq$ the order $m\prec m'\prec\cdots\prec 2\prec2'\prec1\prec 1'$, see also \cite{Oe2}, then a subset of $\un n$ is symplectic standard in the sense of King with the ordering $\preceq$ if and only if it is admissible in the sense of \cite{DeC} (with the natural ordering), see \cite{Sh}.
If $X_R$ is an affine scheme over $R$, $M$ and $N$ $R$-modules and $\mu:X_R\times M\to N$ a morphism such that $\mu(x,\cdot):M(S)\to N(S)$ is linear for all $x\in X(S)$ and all $R$ algebras $S$. Then we get an $R$-linear map $\Delta:M\to N\ot R[X_R]$ by the same recipe as in \cite[I.2.8]{Jan}. If $M$ and $N$ are $R$-algebras and the maps $\mu(x,\cdot)$ are algebra homomorphisms, then $\Delta$ is a homomorphism of algebras. There is of course also a version for a ``right action". In the latter case we can take $X$ to be the $\mb Z$-scheme of $2r\times m$-matrices with totally singular column space from \cite{DeC}. Then we can take $M=\bigwedge\mb Z^{2r}$ and $N=\bigwedge\mb Z^m$. Next we show that the ``action'' $\mu$ kills the $z_t\in \bigwedge\mb Z^{2r}$, $t>0$. Then one applies the comodule map $\Delta$ to \eqref{eq.v} and obtains the identities needed for the straightening and one obtains the spanning results for Kings standard tableaux. Next one proves independence over $\mb Z$ using Lemma~\ref{lem.dim} and the presence of highest weight vectors and finally one proves reducedness as in for $k[M]/(d)$ in Proposition~\ref{prop.CM}(i) in the next section.
%
\end{rem}

\section{A bideterminant basis for a reductive monoid}\label{s.bideterminant basis}
Let $R$ be a commutative ring. Let $M(R)$ be the set of $n\times n$ matrices
$
\begin{bmatrix}
A&0\\
0&B
\end{bmatrix}
$
with $A,B\in\Mat_m(R)$ and $A^TB=AB^T=d(A,B)I$ for some scalar $d(A,B)\in R$. Here $A^T$ denotes the transpose of a matrix $A$ and $I$ is the $m\times m$ identity matrix.
Clearly $M(R)$ is a submonoid of $\Mat_m(R)\times\Mat_m(R)\subseteq\Mat_n(R)$. We denote the group of invertible elements of $M(R)$ by $\tilde G(R)$. Note that $\tilde G(R)$ contains the group $G(R)$ of matrices
$
\begin{bmatrix}
A&0\\
0&(A^{-1})^T
\end{bmatrix},
$
$A\in\GL_m(R)$. We denote the functors $R\mapsto M(R)$, $R\mapsto \tilde G(R)$ and $R\mapsto G(R)$ by $M_\mb Z$, $\tilde G_\mb Z$ and $G_\mb Z$. The functors $M_\mb Z$ and $G_\mb Z$ are closed subschemes of $\Mat_{m,\mb Z}\times\Mat_{m,\mb Z}$ and the functor $\tilde G_\mb Z$ is a closed subscheme of $\GL_{m,\mb Z}\times\GL_{m,\mb Z}$. From now on we will take $J$ to be the second matrix in \eqref{eq.J} and we write $\mb Z^n=V_{\mb Z}\oplus W_{\mb Z}$, where $V_{\mb Z}$ is the sub $\mb Z$-modules spanned by the natural basis elements $v_i$ with $i\in\un m$ and $W_{\mb Z}$ is spanned by the elements $w_i:=v_{i'}$, $i\in\un m$. Then $M_{\mb Z}$ is a closed submonoid scheme of $\SpM_{n,\mb Z}$. Note furthermore that $V_{\mb Z}$ and $W_{\mb Z}$ are sub $\Mat_{m,\mb Z}\times\Mat_{m,\mb Z}$-modules and therefore also sub $M_{\mb Z}$-modules of $\mb Z^n$. For $i,j\in\un m$, we define $h_{ij}$ and $\ov h_{ij}$ to be the restrictions of $g_{ij'}$ and $\ov g_{ij'}$, respectively, to $\Mat_{m,\mb Z}\times\Mat_{m,\mb Z}$. For $i,j\in\un m$ we denote the restriction of $x_{ij}$ to $\Mat_{m,\mb Z}\times\Mat_{m,\mb Z}$ again by $x_{ij}$, but the restriction of $x_{i'j'}$ to $\Mat_{m,\mb Z}\times\Mat_{m,\mb Z}$ will be denoted by $y_{ij}$. So we have for $i,j\in\un m$
\begin{equation*}
h_{ij}=\sum_{l=1}^mx_{li}y_{lj}\text{\quad and\quad}\ov{h}_{ij}=\sum_{l=1}^mx_{il}y_{jl}\ .
\end{equation*}
The ideal of $G_\mb Z$ in $\mb Z[\Mat_{m,\mb Z}\times\Mat_{m,\mb Z}]$ is generated by the elements $h_{ij}$, $1\le i,j\le m$, $i\ne j$, and $h_{rr}-1$, $1\le r\le m$.
The ideal of $\tilde G_\mb Z$ in $\mb Z[\GL_{m,\mb Z}\times\GL_{m,\mb Z}]$ is generated by the elements $h_{ij}$, $1\le i,j\le m$, $i\ne j$, and $h_{rr}-h_{ss}$, $1\le r,s\le m$.
The ideal of $M_\mb Z$ in $\mb Z[\Mat_{m,\mb Z}\times\Mat_{m,\mb Z}]$ is generated by the elements
\begin{equation}\label{eq.M}
\{h_{ij},\,\ov h_{ij},\,h_{rr}-\ov h_{ss}\,|\,1\le i,j\le m, i\ne j, 1\le r,s\le m\}.
\end{equation}
The algebra $\mb Z[\Mat_{m,\mb Z}\times\Mat_{m,\mb Z}]$ is $\mb Z\times\mb Z$-graded by ${\rm deg}(x_{ij})=(1,0)$ and ${\rm deg}(y_{ij})=(0,1)$. The algebra $\mb Z[M_\mb Z]$ inherits a grading, since the ideal of $M_\mb Z$ is generated by homogeneous elements.
The restriction of the coefficient of dilation $d\in\mb Z[\SpM_{n,\mb Z}]$ to $M_\mb Z$ coincides (of course) with the function $d\in\mb Z[M_\mb Z]$ defined above. Note that $G_\mb Z\cong \GL_{m,\mb Z}$. Precisely as in the case of $\GSp_{n,\mb Z}$ in Section~\ref{s.symplecticstraightening} one shows that the fibers of $\tilde G_{\mb Z}$ over $\Spec(\mb Z)$ are reduced. So $\tilde G(K)$ is a closed subgroup of $\GL_n(K)$, defined over the prime field. Furthermore, it is connected and reductive, and each of its elements is a scalar multiple of an element in $G(K)$. The group $\tilde G=\tilde G(k)$ is the group of $k$-points of $\tilde G(K)$ and it is dense in $\tilde G(K)$. It is our aim in this section to show that the functions \eqref{eq.M} generate the vanishing ideal of $M$ in $k[\Mat_m\times\Mat_m]$. This will be deduced from a bitableaux basis result.

The group scheme $G_{\mb Z}$ has the maximal torus of diagonal matrices in common with $\Sp_{n,\mb Z}$, we denote its character group by $X$, it is isomorphic to $\mb Z^m$. We will use the same notation for the restriction of characters of diagonal matrices in $\GL_{n,\mb Z}$ to those in $G_{\mb Z}$ as in Section~\ref{s.symplecticstraightening}. Furthermore, we will again embed $\mb Z^m$ into $\mb Z^n$ by extending with zeros. We will now use the root system of $A_{m-1}$ and our choice of positive roots is the usual one: $\ve_i-\ve_j$, $1\le i<j\le m$. The corresponding Borel subgroup of $G_{\mb Z}$ is the subgroup of matrices of the form
$
\begin{bmatrix}
A&0\\
0&B
\end{bmatrix}
$, with $A$ upper triangular and $B$ lower triangular. The type $A_{m-1}$ order or dominance order on weights and the set of (polynomial) dominant weights are defined as in Section~\ref{s.prelim} with $n$ replaced by $m$. We denote the length of a partition $\lambda$ by $l(\lambda)$.

We assume given a linear order $\preceq$ on $\un m$. Define $\zeta:\un m\to\un m$ by $$\zeta(j)=|\{i\in\un m\,|\,i\preceq j\}|.$$ Note that $i\preceq j$ if and only if $\zeta(i)\le\zeta(j)$. Let $\lambda^1$ and $\lambda^2$ be partitions of length $\le m$. Following Stembridge \cite{St} we define a {\it rational tableau of shape $(\lambda^1,\lambda^2)$} to be a pair $T=(T^1,T^2)$ where $T^i$ is a tableau of shape $\lambda^i$ with entries in $\un m$. For $j\in\un m$, we define the $j^{\rm th}$ column of $T$ to be the pair $(C^1,C^2)$, where $C^i$ is the $j^{\rm th}$ column of $T^i$. Of course $C^1$ or $C^2$ or both may be empty. A pair of subsets $I=(I^1,I^2)$ is called {\it standard} if
$$|\{i\in I^1\,|\,i\preceq j\}|+|\{i\in I^2\,|\,i\preceq j\}|\le\zeta(j)$$
for all $j\in\un m$. We identify each pair of subsets $I=(I^1,I^2)$ of $\un m$ with the one column rational tableau $(T^1,T^2)$ such that, for $i=1,2$ the entries of $T^i$ are the elements of $I^i$ and the entries of $T^i$ are strictly increasing (according to $\preceq$) from top to bottom. A rational tableau is called standard if $T^1$ and $T^2$ are standard and if its first column is standard. If $T$ is standard, then every column of $T$ is standard. Note that if a rational tableau of shape $(\lambda^1,\lambda^2)$ is standard, then $l(\lambda^1)+l(\lambda^2)\le m$. For $\lambda^1, \lambda^2\in\mb Z^m$ we put $$[\lambda^1,\lambda^2]:=\lambda^1-\lambda^{2\,\rm rev}\, ,$$ where $\lambda^{2\,\rm rev}$ is the reversed tuple of $\lambda^2$. It is easy to see that for any $\lambda\in X^+$ there exists unique partitions $\lambda^1$ and $\lambda^2$ with $l(\lambda^1)+l(\lambda^2)\le m$ and $\lambda=[\lambda^1,\lambda^2]$. In the sequel, when $\lambda^1$ and $\lambda^2$ are introduce after $\lambda$, they are supposed to have these properties. If $\lambda\in X^+$ and $\lambda=[\lambda^1,\lambda^2]$ as above, then we say that a rational tableau has shape $\lambda$ if it has shape $(\lambda^1,\lambda^2)$. Let $T=(T^1,T^2)$ be a rational tableau of shape $(\lambda^1,\lambda^2)$, then we define the {\it weight} $\mu$ of $T$ by $\mu=\mu^1-\mu^2$, where $\mu^i$ is the weight of $T^i$ as defined in Section~\ref{s.prelim} (with $n$ replaced by $m$). If the entries in each column of $T^1$ and $T^2$ are distinct, then we have $-[\lambda^2,\lambda^1]=[\lambda^1,\lambda^2]^{\rm rev}\le\mu\le [\lambda^1,\lambda^2]$. 
Now let $\lambda\in X^+$ and write $\lambda=[\lambda^1,\lambda^2]$ as above. Then we define the {\it canonical rational tableau} $T_\lambda$ as the rational tableau $(T^1,T^2)$ of shape $\lambda$ such that $T^1$ has all its entries in the $i^{\rm th}$ row equal to $i$, $1\le i\le l(\lambda^1)$, and $T^2$ has all its entries in the $i^{\rm th}$ row equal to $m-i+1$, $1\le i\le l(\lambda^2)$. Note that $T_\lambda$ has weight $\lambda$ and that one can make it standard for the natural order by reversing each column of $T^2$.

Let $\lambda^1$ and $\lambda^2$ be partitions of length $\le m$. A {\it rational bitableau} of shape $(\lambda^1,\lambda^2)$ is a pair $(S,T)$ where $S$ and $T$ are rational tableaux of shape $(\lambda^1,\lambda^2)$, we call it standard if both $S$ and $T$ are standard. Now let $(S,T)$, $S=(S^1,S^2)$, $T=(T^1,T^2)$, be a rational bitableau of shape $(\lambda^1,\lambda^2)$. Then we define the bideterminant $(S\,|\,T)$ associated to $(S,T)$ by $$(S\,|\,T)=(S^1\,|\,T^1)_1(S^2\,|\,T^2)_2\ ,$$ where $(S^1\,|\,T^1)_1$ is defined by \eqref{eq.bideterminant} with $S$, $T$ and $\lambda$ replaced by $S^1$, $T^1$ and $\lambda^1$ and $(S^2\,|\,T^2)_2$ is defined by \eqref{eq.bideterminant} with $x_{ij}$, $S$, $T$ and $\lambda$ replaced by $y_{ij}$, $S^2$, $T^2$ and $\lambda^2$. We define the {\it weight} of a rational bitableau $(S,T)$ to be $(-\mu,\nu)\in X\times X$, where $\mu$ is the weight of $S$ and $\nu$ is the weight of $T$. If $H_{\mb Z}$ is the maximal torus of diagonal matrices of $G_\mb Z$, then the bideterminant $(S\,|\,T)$ is an $H_{\mb Z}\times H_{\mb Z}$ weight vector with weight equal to that of $(S,T)$. The {\it degree} of a rational tableau or bitableau of shape $(\lambda^1,\lambda^2)$ is defined to be $(|\lambda^1|,|\lambda^2|)$. Note that the degree of a bideterminant $(S\,|\,T)$ is equal to that of the rational bitableau $(S,T)$. The {\it content} of a rational tableau $T=(T^1,T^2)$ is the tuple $a\in\mb Z^m$, such that $a_i$ is the number of occurrences of $i$ in $T^1$ and $T^2$.

We give the exterior algebra $\bigwedge(V_{\mb Z}\oplus W_{\mb Z})$ a $\mb Z\times\mb Z$-grading by giving the elements $v_i$ degree $(1,0)$ and the elements $w_i$ degree $(0,1)$. For $I=(I^1,I^2)\in P(m,r)\times P(m,s)$ we define $v_I:=v_{I^1}\wedge w_{I^2}\in\bigwedge(V_{\mb Z}\oplus W_{\mb Z})^{r,s}$, where $v_{I^1}$ and $w_{I^2}$ are defined as in Section~\ref{s.symplecticstraightening} (for $w_{I^2}$ we use the $w_i$). Note that $v_I$ as defined here is equal to $\pm v_{I^1\cup (I^2)'}$ as defined in Section~\ref{s.symplecticstraightening}, where $(I^2)'=\{i'\,|\,i\in I^2\}$. In particular, the $v_I$, $I\in P(m,r)\times P(m,s)$, form a basis of $\bigwedge(V_{\mb Z}\oplus W_{\mb Z})^{r,s}$. The elements $z_r$ and the order $\preceq$ on $\mb Z^m$ is defined as in Section~\ref{s.symplecticstraightening}. We observe that if $I,J\in P(n,r)$ have the same degree ($|I|=|J|$) and the same symplectic weight, then $|I\cap\un m|=|J\cap\un m|$ and $|I\sm\un m|=|J\sm\un m|$. It is now clear that one can reformulate Proposition~\ref{prop.extbasis} with $P(n,r)$ replaced by $P(m,r)\times P(m,s)$ and $P(n,r-t)$ replaced by $P(m,r-t)\times P(m,s-t)$. The action of $\Mat_{m,\mb Z}\times\Mat_{m,\mb Z}$ on $\bigwedge(V_{\mb Z}\oplus W_{\mb Z})$ stabilises the $\mb Z\times\mb Z$-grading and for the comodule map $\Delta_\wedge:\bigwedge(V_{\mb Z}\oplus W_{\mb Z})^{r,s}\to\bigwedge(V_{\mb Z}\oplus W_{\mb Z})^{r,s}\ot\mb Z[M_\mb Z]$ of the $M_\mb Z$-action we have
\begin{equation}\label{eq.comod2}
\Delta_\wedge(v_J)=\sum_Iv_I\ot(I\,|\,J)\, ,
\end{equation}
where the sum is over all $I\in P(m,r)\times P(m,s)$. This just follows by restriction from the corresponding equations for the comodule map of the $\Mat_{m,\mb Z}\times\Mat_{m,\mb Z}$-action which in turn follows from our new definition of the $v_I$ and of the bideterminants and from the fact that the comodule map is a homomorphism of algebras. Next we note that Lemma~\ref{lem.semiinv} also holds for the $M_\mb Z$-action, since $M_\mb Z$ is a closed subscheme of $\SpM_{n,\mb Z}$. Define the integers $c_{LI}$, $I\in P(m,r)\times P(m,s)$, $L\in P(m,r-t)\times P(m,s-t)$,  by $z_t\wedge v_L=\sum_Ic_{LI}v_I$. Now we obtain the following corollary of which the proof is completely analogous to that of Corollary~1 to Proposition~\ref{prop.extbasis}. To obtain \eqref{eq.tab4} from \eqref{eq.tab3} we use the automorphism which sends $x_{ij}$ to $x_{ji}$ and $y_{ij}$ to $y_{ji}$.
\begin{mycor}
Let $I,J\in P(m,r)\times P(m,s)$ with $J$ not standard. Then we have in $\mb Z[M_\mb Z]$
\begin{align}
(I\,|\,J)=\sum_L a_{JL}(I\,|\,L)+\sum_{t,L,L'} b_{J,L}c_{L'I}d^t(L'\,|\,L)&\text{\quad and}\label{eq.tab3}\\
(J\,|\,I)=\sum_L a_{JL}(L\,|\,I)+\sum_{t,L,L'} b_{J,L}c_{L'I}d^t(L\,|\,L')&\label{eq.tab4},
\end{align}
with $a_{JL},b_{JL}\in\mb Z$ given by \eqref{eq.v}; in both cases the first is sum over all $L\in P(m,r)\times P(m,s)$ standard
with $L\vartriangleright J$ and the second sum is over all $t\in\{1,\ldots,\lfloor r/2\rfloor\}$ and $L,L'\in P(m,r-t)\times P(m,s-t)$ standard. Furthermore, all the $L$ and $L'$ occurring have the same weight as $J$ and $I$ respectively.
\end{mycor}

\begin{rems}\label{rem.weights}
1.\ In our new labeling for the standard basis elements of $\bigwedge\mb Z^n$ one can formulate a stronger version of Proposition~\ref{prop.extbasis} and the above corollary. This is based on \cite[Cor~1.2, Prop.~1.3]{DeCS}. Instead of using the identity of Remark~\ref{rems.sympstraight}.3 one uses the identity $$\sum_{L\in P(m,r), L\cap K=\emptyset}z(L)\equiv0\quad({\rm mod}\, N)$$
for $K\subseteq\un{m}$ with $|K|<r$. This is done as follows. Define the ordering $\preceq$ on $P(m,r)$ by $I=\{i_1\prec\cdots\prec i_r\}\preceq \{j_1\prec\cdots\prec j_s\}=J$ if $r\ge s$ and $i_t\preceq j_t$ for all $t\in\{1,\ldots,s\}$. This is equivalent to $\nu_i(I)\ge\nu_i(J)$ for all $i\in\un{m}$, where $\nu_i(I)=|\{j\in I\,|\,j\preceq i\}|$. See e.g. \cite{St}. Now assume $I=(I^1,I^2)\in P(m,r)\times P(m,s)$ is not standard. Pick $i$ minimal with $\nu_i(I^1)+\nu_i(I^2)>i$. Then $i\in I^1\cap I^2$ and $\nu_{i-1}(I^1)+\nu_{i-1}(I^2)=i-1$, i.e. $\nu_{i-1}(I^1)=\nu_{i-1}((I^2)^c)$, where $(I^2)^c$ is the complement of $I^2$ in $\un{m}$. Now put $J=I^1\cap I^2\cap[1,i]$ and $K=(I^1)^c\cap (I^2)^c\cap[1,i)$, where we use the standard interval notation (for our ordering $\preceq$ on $\un{m}$). Then $J\cap K=\emptyset$ and one easily deduces that $r:=|J|=|K|+1$. 
Now $v_I=\pm v_{I'}z(J)$, where $I'=(I^1\sm J,I^2\sm J)$, so we can apply the above identity to the factor $z(J)$ of $v_I$. Then we may restrict the summation to those $L$ that are also disjoint from $I^1\sm J$ and $I^2\sm J$, i.e. to those $L$ that are subsets of $J\cup\big((I^1\cup I^2)^c\cap(i,\infty)\big)$. The final result is that Proposition~\ref{prop.extbasis} with the new labelling and the above corollary also hold if we replace ``$L\vartriangleright J$" by ``$L^1\succ J^1$ and $L^2\succ J^2$".\\
2.\ Assume $k=K=\mb C$ and let $\Lambda_{r,s}\subseteq\mb Z^m$ be the set of weights $[\lambda^1,\lambda^2]$ where $\lambda^1$ and $\lambda^2$ are partitions with $l(\lambda^1)+l(\lambda^2)\le m$, $|\lambda^1|\le r$, $|\lambda^2|\le s$ and $r-|\lambda^1|=s-|\lambda^2|$. Then one easily checks that these are the dominant weights of the $\GL_n$-module $V^{\ot r}\ot(V^*)^{\ot s}$. Now one can show that each of these weights is actually the weight of a highest weight vector. Using the fact that $V\ot V^*$ contains the trivial module one is reduced to the case that $r-|\lambda^1|=s-|\lambda^2|=0$ and then one can simply tensor two highest weight vectors together. This is standard, see e.g. \cite{St} and \cite{BCHLLS}, where also multiplicity questions are considered. As is well known, the dominant weights of the irreducible $\GL_n$-module with highest weight $\lambda$ are the dominant weights $\le\lambda$, so it follows from the above remarks that $\Lambda_{r,s}$ is saturated.
\end{rems}

To prove the Theorem~\ref{thm.M} we need the following notation. For $\lambda,\mu\in\mb Z^m$ we write $\mu\le_1\lambda$ if $\sum_{i=1}^j\mu_i\le\sum_{i=1}^j\lambda_i$ for all $j\in\un m$.

\begin{lem}\label{lem.part}
Let $\mu^1,\mu^2\in\mb Z^m$ have entries $\ge0$, let $\lambda^1$ and $\lambda^2$  be partitions with $l(\lambda^1)+l(\lambda^2)\le m$ and assume that $\mu^i\le_1\lambda^i$ for $i=1,2$ and that $|\lambda^1|-|\mu^1|=|\lambda^2|-|\mu^2|$. Then $[\mu^1,\mu^2]\le[\lambda^1,\lambda^2]$.
\end{lem}

\begin{proof}
Put $t=|\lambda^1|-|\mu^1|=|\lambda^2|-|\mu^2|$ and, for $i=1,2$, $l_i=l(\lambda^i)$, $\nu^i=\lambda^i-t\e_{l_i}$. Then we have for $i=1,2$, $\mu^i\le_1\nu^i$
and $|\mu^i|=|\nu^i|$, so $\mu^i\le\nu^i$. Now $\eta\mapsto\eta^{\rm rev}$ reverses the order $\le$ (but {\it not} $\le_1$) and $\eta\mapsto-\eta$ also reverses this order, so $[\mu^1,\mu^2]\le[\nu^1,\nu^2]$. But $[\nu^1,\nu^2]=[\lambda^1,\lambda^2]-t(\e_{l_1}-\e_{m-l_2+1})\le[\lambda^1,\lambda^2]$, since $(e_j)^{\rm rev}=e_{m-j+1}$ and $l_1<m-l_2+1$.
\end{proof}

\begin{thm}\label{thm.M}
Recall that the type $A_{m-1}$ partial order on weights is denoted by~$\le$.
\begin{enumerate}[{\rm (i)}]
\item Let $\lambda\in X^+$ be a dominant weight and let $S$ and $T$ be rational tableaux of shape $\lambda$. Then the bideterminant $(S\,|\,T)\in \mb Z[M_\mb Z]$ can be written as a linear combination $\sum_ia_id^{t_i}(S_i\,|\,T_i)$, where the $a_i$ are integers and $(S_i,T_i)$ is standard of shape $\mu_i\le\lambda$ with $(|\mu_i^1|,|\mu_i^2|)=(|\lambda^1|,|\lambda^2|)-(t_i,t_i)$ and with the same weight as $(S,T)$.
\item The elements $d^t(S\,|\,T)$ with $t$ an integer $\ge0$ and $S$ and $T$ standard form a basis of $\mb Z[M_\mb Z]$.
\end{enumerate}
\end{thm}

\begin{proof}
(i).\ By Lemma~\ref{lem.part} it is enough to prove that for $\lambda^1$ and $\lambda^2$ partitions of length $\le m$ and $S$ and $T$ rational tableaux of shape $(\lambda^1,\lambda^2)$ the bideterminant $(S\,|\,T)$ can be written as $\sum_ia_id^{t_i}(S_i\,|\,T_i)$, where the $a_i$ are integers and $(S_i,T_i)$ is standard of shape $(\mu_i^1,\mu_i^2)$ with $\mu_i^1\le_1\lambda^1$, $\mu_i^2\le_1\lambda^2$ and $(|\mu_i^1|,|\mu_i^2|)=(|\lambda^1|,|\lambda^2|)-(t_i,t_i)$. First we show that, when both parts of $S$ and $T$ are $\GL_m$-standard, but $S$ or $T$ is not standard, $(S\,|\,T)$ can be written as a $\mb Z$-linear combination $\sum_ia_i(S_i\,|\,T_i)+\sum_ib_id^{t_i}(S'_i\,|\,T'_i)$, where in the first sum $S_i$ and $T_i$ are of shape $(\lambda^1,\lambda^2)$ and $S_i$ or $T_i$ has content $\vartriangleright$ than that of $S$ or $T$ respectively and in the second sum the $t_i$ are $>1$ and the $S'_i$ and $T'_i$ are of shape $(\mu_i^1,\mu_i^2)$ with $\mu_i^1\le_1\lambda^1$, $\mu_i^2\le_1\lambda^2$ and $(|\mu_i^1|,|\mu_i^2|)=(|\lambda^1|,|\lambda^2|)-(t_i,t_i)$.
To prove this we apply \eqref{eq.tab3} or \eqref{eq.tab4} to the product of two minors corresponding to the first columns of $S$ and $T$, and all that remains to check is that for the new shapes $(\mu^1,\mu^2)$ that arise from the second sums in \eqref{eq.tab3} or \eqref{eq.tab4} we have  $\mu^1\le_1\lambda^1$, $\mu^2\le_1\lambda^2$ and $|\lambda^1|-|\mu^1|=|\lambda^2|-|\mu^2|$. This follows easily from the fact that such a shape $\mu$ is formed by shortening the length of the first columns of $\lambda^1$ and $\lambda^2$ by the same number and moving the resulting columns to the right positions to get a pair of Young diagrams. Now we can finish by induction as in the proof of Theorem~\ref{thm.SpM}(i) applying the $\GL_m$-straightening separately to the two bitableaux of a rational bitableau.\\
(ii).\ We proceed as in the proof of Theorem~\ref{thm.SpM}(ii). The piece of degree $(r,s)$ of $\mb C[M_\mb C]$ surjects onto that of $\mb C[\ov{\tilde G(\mb C)}]$. By Proposition~\ref{prop.envalg} this is the dual space of the enveloping algebra $E$ of $G(\mb C)$ (or $\tilde G(\mb C)$) in $\End_{\mb C}(V_\mb C^{\ot r}\ot W_\mb C^{\ot s})$. Clearly, $E$ is also the enveloping algebra of $\GL_m(\mb C)$ in $\End_{\mb C}(V_\mb C^{\ot r}\ot(V^*_\mb C)^{\ot s})$. By Weyl's complete reducibility theorem for complex semisimple Lie algebras $E$ has dimension $\sum_\lambda \dim L(\lambda)^2$, where $L(\lambda)$ is the irreducible $\GL_m(\mb C)$-module of highest weight $\lambda$ and the sum is over all $\lambda$ such that $L(\lambda)$ appears in $V_\mb C^{\ot r}\ot(V^*_\mb C)^{\ot s}$. By Remark~\ref{rem.weights}.2 these are the weights in $\Lambda_{r,s}$. Now $\dim L(\lambda)$ is the number of standard rational tableaux of shape $\lambda$ by \cite[Prop.~2.4(a)]{St}, so $\dim E$ is the number of rational bitableaux $(S,T)$ with $S$ and $T$ standard with shape $\in\Lambda_{r,s}$, i.e. with degree $(r-t,s-t)$, $t\in\{0,\ldots,\min(r,s)\}$. So, by (i) and Lemma~\ref{lem.dim}, the elements $d^t(S\,|\,T)$ with $t\in\{0,\ldots,\min(r,s)\}$ and $S$ and $T$ standard of degree $(r-t,s-t)$ form a basis of the degree $(r,s)$ piece of $\mb Z[M_\mb Z]$.
\end{proof}

The arguments for the proof of the corollary below are precisely the same as in the case of $\SpM_{n,\mb Z}$. See the paragraph after the proof of Theorem~\ref{thm.SpM}.

\begin{cornn}\
\begin{enumerate}[\rm(i)]
\item The fibers of $M_\mb Z$ over $\Spec(\mb Z)$ are reduced and irreducible.
\item The group $\tilde G$ is dense in $M(K)$ and the functions \eqref{eq.M} generate the vanishing ideal of $M$ in $k[\Mat_m\times\Mat_m]$.
\end{enumerate}
\end{cornn}

Using the isomorphism $G_{\mb Z}\cong\GL_{m,\mb Z}$ we can consider the $x_{ij}$ and $y_{ij}$ and the bideterminants as functions on $\GL_{m,\mb Z}$. For $\lambda\in X^+$ one defines $A_{\le\lambda},A_{<\lambda},\nabla_{\mb Z}(\lambda),$\\ $\tilde\nabla_{\mb Z}(\lambda)\subseteq\mb Z[\GL_{m,\mb Z}]$ completely analogous to the $\Mat_{n,\mb Z}$-case. Simply replace tableaux by rational tableaux.

\begin{thm}\label{thm.GLm}
The analogue of Theorems~\ref{thm.GLn} and \ref{thm.Sp} holds for $\GL_{m,\mb Z}\times\GL_{m,\mb Z}$ acting on $\mb Z[\GL_{m,\mb Z}]$, arbitrary dominant weights and rational bitableaux.
\end{thm}

\begin{proof}
Since $\mb Z[M_{\mb Z}]/(d-1)\cong G_{\mb Z}\cong\GL_{m,\mb Z}$ (i) and (ii) follow from Theorem~\ref{thm.M}.\\
(iii).\ First we show that the given tableaux span $\nabla_{\mb Z}(\lambda)$ and $\tilde\nabla_{\mb Z}(\lambda)$. This is done precisely as in the proof of Theorem~\ref{thm.Sp}(iii) using the fact that a column $I=(I^1,I^2)$ of $T_\lambda$ always has the property that $I^1\cap I^2=\emptyset$; for such $I$ the second sums in \eqref{eq.tab3} and \eqref{eq.tab4} are zero. To show independence we use the trick of Remark~\ref{rems.GL}.1: multiplying on the left or on the right with a matrix
$\begin{bmatrix}
P&0\\
0&P
\end{bmatrix}$, $P$ a permutation matrix,
we see that (ii) also holds if we use two orderings, $\preceq_1$ for the left tableau and $\preceq_2$ for the right tableau. In case of $\nabla_{\mb Z}(\lambda)$ we can choose $\preceq_1$ such that $T_\lambda$ is standard (since $l(\lambda^1)+l(\lambda^2)\le m$)
and in case of $\tilde\nabla_{\mb Z}(\lambda)$ we can choose $\preceq_2$ such that $T_\lambda$ is standard. Then (ii) gives us independence.
\\
(iv). This is proved precisely as in Theorem~\ref{thm.Sp}(iv).
\end{proof}

Now let $B_{\mb Z}$ and $B_{\mb Z}^-$ be the Borel subgroups of upper resp. lower triangular matrices in $\GL_{m,\mb Z}$. Then the analogue of the corollary to Theorem~\ref{thm.GLn} holds, its proof is also completely analogous.

\begin{cornn}
The analogue of the corollaries to Theorems~\ref{thm.GLn} and \ref{thm.Sp} holds for $\GL_{m,\mb Z}\times\GL_{m,\mb Z}$ acting on $\mb Z[\GL_{m,\mb Z}]$, arbitrary dominant weights and rational bitableaux.
\end{cornn}
Combining restriction of functions with the isomorphism $G_k\cong(\GL_m)_k$ we obtain from Theorem~\ref{thm.M} and the above corollary a canonical isomorphism $O_{\Lambda_{r,s}}(k[\GL_m])\cong k[M]^{r,s}$ of coalgebras and therefore, by Proposition~\ref{prop.envalg}(i), a canonical coalgebra isomorphism of $O_{\Lambda_{r,s}}(k[\GL_m])$ with the dual of the enveloping algebra of $\GL_m$ in $\End_k(V^{\ot r}\ot(V^*)^{\ot s})$. This enveloping algebra is called the {\it rational Schur algebra}, see \cite{DiDo} and \cite{Don5}. From the above remarks it is clear that the rational Schur algebra is a generalised Schur algebra, see \cite[A.16]{Jan}, and therefore quasihereditary.

In the remainder of this section we assume that $k=K$ is algebraically closed. For convenience, we will denote a matrix
$
\begin{bmatrix}
A&0\\
0&B
\end{bmatrix}
\in\Mat_m\times\Mat_m$ by $(A,B)$. Let $X$ be the variety of noninvertible elements in $M$. Clearly, $X$ is the set of zeros of $d$ in $M$, so it consists of the matrices $(A,B)$ where $A,B\in\Mat_m$ are such that $A^TB=AB^T=0$. For $r\in\un m$ we define the idempotents $E_r,F_r\in\Mat_m$ to be the diagonal matrices of which the first resp. last $r$ diagonal entries are equal to $1$ and all other diagonal entries are $0$. Let $r,s\in\un m$ with $r+s\le m$. We define $E_{r,s}\in M$ by
$E_{r,s}=(E_r,F_s)$ and we define $X_{rs}\subseteq X$ to be the variety of matrices $(A,B)\in M$ for which $\rk\,A\le r$ and $\rk\,B\le s$. Here $\rk\,A$ denotes the rank of a matrix $A$. Note that the $\tilde G\times \tilde G$-orbit of an idempotent $E_{r,s}$ is equal to its $G\times G$-orbit. We have $X_{r_1,s_1}\subseteq X_{r_2,s_2}$ if and only if $r_1\le r_2$ and $s_1\le s_2$.

We let $\GL_m\times\GL_m$ act on $M$ using the isomorphism $\GL_m\cong G$. Let $U$ and $U^-$ be the subgroups of $\GL_m$ that consist of, respectively, the upper and lower unitriangular matrices and let $\lambda\in X^+$. Recall that, by the corollary to Theorem~\ref{thm.GLm}, the bideterminants $(T_\lambda\,|\,S)$ are $U^-$-fixed under the left regular action and that the $(S\,|\,T_\lambda)$ and are $U$-fixed under the right regular action. Of course, this also holds when we interpret the bideterminants as functions on $M$. It is a simple exercise in linear algebra to show that every element of $X$ is $\GL_m\times\GL_m$-conjugate to one of the idempotents $E_{r,s}$. The variety $X_{rs}$ is the closure of the $\GL_m\times\GL_m$-orbit of $E_{rs}$ and a simple centraliser computation shows that the dimension of $X_{rs}$ is $m^2-(m-(r+s))^2=(r+s)(2m-(r+s))$.
The varieties $X_{r,m-r}$, $0\le r\le m$, are the irreducible components of $X$. We define $I_{rs}$ to be the ideal of $k[X]$ generated by the minors in the $x_{ij}$ of degree $r+1$ together with the minors in the $y_{ij}$ of degree $s+1$.
\begin{prop}\label{prop.CM}\
\begin{enumerate}[\rm(i)]
\item The element $d$ generates the vanishing ideal of $X$ in $k[M]$ and the bideterminants $(S\,|\,T)$ with $S$ and $T$ standard form a basis of $k[X]$. Furthermore, $M$ is normal and $X$ and $M$ are Cohen-Macaulay.
\item The ideal $I_{r,s}$ is the vanishing ideal of $X_{rs}$ in $k[X]$ and the bideterminants $(S\,|\,T)$ with $S$ and $T$ standard form a basis of $k[X]$ and those whose shape $(\lambda^1,\lambda^2)$ satisfies $l(\lambda^1)\le r$ and $l(\lambda^2)\le s$ form a basis of $k[X_{r,s}]$.
\item(\cite{St}, \cite{MeTr}) The varieties $X_{rs}$ are normal and Cohen-Macauley and have a rational desingularisation.
\end{enumerate}
\end{prop}

\begin{proof}
(i).\ One can define the $\GL_m$ modules $A_{\le\lambda},A_{<\lambda},\nabla(\lambda)$ and $\tilde\nabla(\lambda)$ inside $k[M]/(d)$ and Theorem~\ref{thm.GLm} is then also valid for $k[M]/(d)$. Moreover, $\nabla(\lambda)$ and $\tilde\nabla(\lambda)$ can be defined inside $k[M]$ and then Theorem~\ref{thm.GLm}(iii) is still valid. It now follows that the natural maps $k[M]\to k[\GL_m]$ and $k[M]\to k[M]/(d)$ restrict to isomorphisms between the different versions of $\nabla(\lambda)$ and $\tilde\nabla(\lambda)$. So, by the corollary to Theorem~\ref{thm.GLm}, the versions of these modules inside $k[M]/(d)$ are induced modules. Now we have by \cite[II.2.13,II.4.13]{Jan} for $\lambda,\mu\in X^+$ that
$$
\dim (k[M]/(d))^{U^-\times U}_{\lambda,\mu}=
\begin{cases}
1,\text{\quad if\ } \lambda=-\mu,\\
0,\text{\quad otherwise.}
\end{cases}
$$
Clearly, $(T_\lambda\,|\,T_\lambda)$ is a weight vector of weight $(-\lambda,\lambda)$ and it is $U^-\times U$-fixed. So the vectors $(T_\lambda\,|\,T_\lambda)$, $\lambda\in X^+$ form a basis of $(k[M]/(d))^{U^-\times U}$, each with a distinct weight. So if a $B^-\times B$-submodule of $k[M]/(d)$ is nonzero, then it must contain one of these vectors. In particular this applies to the radical of the algebra $k[M]/(d)$. However, $(T_\lambda\,|\,T_\lambda)$ is nonzero as a function on $X$, since it is nonzero on $E_{r,s}$, where $r=|\lambda^1|$, $s=|\lambda^2|$, $\lambda=[\lambda^1,\lambda^2]$. So $k[M]/(d)$ is reduced and $d$ generates the vanishing ideal of $X$ in $k[M]$. The second assertion is now also clear. We now show that $M$ is normal. Since $k[M][d^{-1}]=k[\tilde G]$ is integrally closed, we only have to show that $k[M]$ is integrally closed in $k[M][d^{-1}]$. So, let $f\in k[M][d^{-1}]$ be integral over $k[M]$. Write $f=f_1/d^r$, with $f_1\in k[M]$ such that $d\nmid f_1$ and assume that $r\ge1$. Now $f$ satisfies some monic equation, of degree $s$ say, with coefficients in $k[M]$. Multiplying through with $d^{rs}$ we get that $d|f_1^s$. Since, by (ii), the ideal $(d)$ is radical, this implies that $d|f_1$. This contradicts our choice of $f_1$, so $r=0$ and $f=f_1\in k[M]$. Now $M$ is Cohen-Macaulay by \cite[Cor~2]{Rit} or \cite[Cor.~6.2.9]{BriKu} and $X$ is then also Cohen-Macaulay, since $k[X]=k[M]/(d)$.\\
(ii).\ We first work over $\mb Z$. Clearly, $\mb Z[X]/I_{rs,\mb Z}$ is spanned by the given bideterminants. Then we deduce from Lemma~\ref{lem.dim} and the presence of the highest weight vectors $(T_\lambda\,|\,T_\lambda)$ that they must form a basis. This holds then also over $k$ by base change and we can show that $k[X]/I_{rs}$ is reduced as in the case of $k[X]/(d)$ in (ii).\\
(iii).\ All assertions are proved in \cite{MeTr} using Frobenius splitting and Schubert varieties (one has to use the isomorphism $(A,B)\mapsto(A^T,B)$ between $X$ and the variety in \cite{MeTr} or \cite{Str}). The first two assertions are proved in \cite{Str} using Hodge algebras (see e.g. \cite[Ch.~7]{BrHe} for a definition). The problem is however that the partial order on the minors given there does not satisfy axiom (H2), see Remark~\ref{rems.hodge}.1. We indicate how one can repair this using Gr\"obner bases. In the polynomial algebra on the minors on the first and second matrix of size $\le r$ and $\le s$ respectively, consider the ideal generated by the $\GL_m$-straightening relations for the minors on both matrices separately (the minors in $I_{rs}$ are understood to be zero) and by the relations \eqref{eq.tab3},\eqref{eq.tab4} (note that in these equations a symbol $(I\,|\,J)$ denotes a product of two minors and that the second sums are now zero). Now we order the minors according to the ordering $\succeq$ from Remark~\ref{rems.sympstraight}.1, keeping minors on different matrices incomparable. Then we extend this order to a linear order. Now we choose the reverse lexicographic order based on this order of the variables as our monomial order. Then the above relations form a Gr\"obner basis. This follows from Remark~\ref{rem.weights}.1, \cite[Thm.~14.6]{MiSt} and the fact that if a term on the RHS of a $\GL_m$-straightening relation does not satisfy (H2) for our order $\succeq$, then it must be a single minor. The latter follows from \cite[Lemma~2.2]{DeCEP}, since the homogeneous coordinate ring of the Grassmannian is also a Hodge algebra for the opposite of the usual order on the minors. Now we grade the algebra by the size of the minors and we apply \cite[Cor.~8.31]{MiSt} and deduce that $k[X_{rs}]$ is Cohen-Macauley, since this holds by \cite[Prop.~2.6]{Str} for the discrete algebra. Normality can now also be deduced as in \cite{Str}.
\end{proof}

\begin{rems}\label{rems.hodge}
1.~The notion of standardness in \cite{DeCS} and \cite{Str} is equivalent to ours. In \cite[Def.~2.4]{Str} the ordering of the minors is such that on the minors of the second matrix it is the reverse of the usual one. However, with this order $k[X]$ is not a Hodge algebra. This can already be seen by taking $m=2$ and considering the straightening relation $y_{12}y_{21}=y_{11}y_{22}-\det$. Then $\det\le y_{12},y_{21}$ in the usual order, so it cannot also have this property for the opposite order. In \cite{DeCS} a partial order on the minors for which (H2) should hold is not clearly specified; one can repair things in the same way as in the proof of Proposition~\ref{prop.CM}. Finally, we note that the minors in \cite[III.16]{DeCEP} do not form a basis. The notion of standardness there amounts to the following alternative definition of rational standard tableau: A tableau $(S^1\,|\,S^2)$ is standard if $S^1$ and $S^2$ are standard in the usual way and for each column $C^1$ of $S^1$ and each column $C^2$ of $S^2$ one has $C^1\sqsubseteq C^2$ in the ordering on sets based on the opposite ordering of the natural ordering of $\un{m}$. With this notion of standardness we get, with $m=4$, $9$ standard rational tableaux of shape $((2,1),(2,1))$ and weight $(1,0,0,-1)$, but there should only be $5$.\\
2.~Clearly the element $d\in k[M]$ is irreducible if $m\ge2$, so $k[M]$ is not a unique factorisation domain if $m\ge2$.\\
3.~Similar bideterminants as the ones in this section have been considered in \cite{DiDoSt}. It is not clear to me whether the statements in this section about the new shapes that show up during straightening hold in the quantum setting. In view of \cite{Oe2} it seems plausible that there exists a quantum version of our monoid $M$ whose coordinate ring (that is of course the real quantum monoid) has a bideterminant basis involving powers of a quantum coefficient of dilation. This should give an alternative proof of the double centraliser theorem in \cite{DiDoSt}.
\end{rems}


\section{The double centraliser theorem for the rational Schur algebra and the walled Brauer algebra}\label{s.doublecentraliser}
Let $r,s$ be integers $\ge0$. For any $\delta\in k$ one has the Brauer algebra $B_{r+s}(\delta)$ over $k$, see \cite{Br} or one of the many papers in the literature on Brauer algebras. This also makes sense for $\delta$ an integer, since we can replace that integer by its natural image in $k$. A {\it walled Brauer diagram} is a Brauer diagram in which the vertical edges join one of the first $r$ vertices in the top row with one of the first $r$ vertices in the bottom row or one of one of the last $s$ vertices in the top row with one of the last $s$ vertices in the bottom row, and in each row the horizontal edges join one of the first $r$ vertices with one of the last $s$ vertices. So if we draw a wall after the first $r$ vertices in both rows as follows
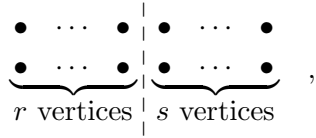
\vspace{.3cm}
\begin{equation*}
\begin{xy}
(-8.4,2.8)*{%
\xymatrix @R=14pt @C=14pt @M=-2pt{
{\bullet}&\cdots&{\bullet}&{\bullet}&\cdots&{\bullet}\\
{\bullet}&\cdots&{\bullet}&{\bullet}&\cdots&{\bullet}
}},
(-9.9,2.8)*=<47pt,28pt>{}*\frm{_\}},%
(-9.9,-5.2)*{\text{$r$ vertices}},
(8.8,2.8)*=<47pt,28pt>{}*\frm{_\}},%
(8.8,-5.2)*{\text{$s$ vertices}},
(-0.3,4.5)*{%
\xymatrix @R=50pt @C=14pt @M=-2pt{
{}\ar@{--}[1,0]\\
{}
}}
\end{xy}\quad,
\end{equation*}
\vspace{.2cm}

\noindent then the horizontal edges must cross the wall and the vertical edges must stay on one side of the wall.
The {\it walled Brauer algebra} $B_{r,s}(\delta)$, see \cite{BCHLLS} or \cite{CdVMD}, is defined as the span in $B_{r+s}(\delta)$ of the walled Brauer diagrams. It is a simple matter to check that this is indeed a subalgebra of $B_{r+s}(\delta)$.

The standard basis of $V=k^m$ determines standard bases of $V^{\ot r}$ and of $V^{\ot r}\ot V^{\ot s}$. We denote the entry of index $((i_1,\ldots,i_r),(j_1,\ldots,j_r))$ of the matrix of the endomorphism of $V^{\ot r}$ given by $A\in\Mat_m$ with respect to the standard basis by $a_{i_1\cdots i_r,j_1\cdots j_r}$. Then the entry of index $$\big(((i_1,\ldots,i_r),(u_1,\ldots,u_s)),((j_1,\ldots,j_r),(v_1,\ldots,v_s))\big)$$ of the matrix of the endomorphism of $V^{\ot r}\ot V^{\ot s}$ given by $(A,B)\in\Mat_m\times\Mat_m$ with respect to the standard basis equals
$a_{i_1\cdots i_r,j_1\cdots j_r}b_{u_1\cdots u_s,v_1\cdots v_s}$. Here $A$ acts on the first $r$ tensor factors and $B$ acts on the last $s$ tensor factors. When we consider $V^{\ot r}\ot V^{\ot s}$ as a $\GL_m$-module via the embedding $\GL_m\hookrightarrow M\subseteq\Mat_m\times\Mat_m$ given in Section~\ref{s.bideterminant basis}, then we write $V^{\ot r}\ot (V^*)^{\ot s}$ instead of $V^{\ot r}\ot V^{\ot s}$.

\begin{prop}\label{prop.envalgeqs}
The enveloping algebra of $\GL_m$ in $\End_k(V^{\ot r}\ot(V^*)^{\ot s})$ is, within the enveloping algebra of $\GL_m\times\GL_m$, defined by the equations
\begin{align}\label{eq.envalg}
\delta_{i_1,u_1}\sum_{l=1}^ma_{li_2\cdots i_r,j_1\cdots j_r}b_{lu_2\cdots u_s,v_1\cdots v_s}=
\delta_{j_1,v_1}\sum_{l=1}^ma_{i_1\cdots i_r,lj_2\cdots j_r}b_{u_1\cdots u_s,lv_2\cdots v_s}\, ,
\end{align}
for $i_1,\ldots,i_r,u_1,\ldots,u_r,j_1,\ldots,j_s,v_1,\ldots,v_s\in\un m$.
\end{prop}

\begin{proof}
By (i) of the corollary to Theorem~\ref{thm.M} we have that the elements \eqref{eq.M}, i.e. the elements $\delta_{i,u}h_{jv}-\delta_{j,v}\ov{h}_{iu}$, generate the vanishing ideal of $M$ in $k[\Mat_m\times\Mat_m]$. Now the corollary to Proposition~\ref{prop.envalg} gives the assertion, since the equation \eqref{eq.envalg} corresponds under the isomorphism $\eta$ from the corollary to the element $\delta_{i_1,u_1}h_{j_1v_1}-\delta_{j_1,v_1}\ov{h}_{i_1u_1}$ multiplied by the monomial $x_{i_2j_2}\cdots x_{i_rj_r}y_{u_2v_2}\cdots y_{u_sv_s}$.
\end{proof}

Now we fix $r+s$ vector symbols $x_1,\ldots,x_{r+s}$ and $r+s$ covector symbols $y_1,\ldots,y_{r+s}$. We consider $y_i$ as the $i^{\rm th}$ component function $x_i$ as the $r+s+i^{\rm th}$ component function on $\bigoplus^{r+s}V^*\oplus\bigoplus^{r+s}V$. We put $\la f,x\ra=\la x,f\ra=f(x)$ for $f\in V^*$ and $x\in V$. With a walled Brauer diagram we associate a $2(r+s)$-multilinear function $F(D)$ on $\bigoplus^{r+s}V^*\oplus\bigoplus^{r+s}V$ as follows. We label the vertices in the top row from left to right with $y_1,\ldots,y_r,x_{r+1},\ldots,x_{r+s}$ and the vertices in the bottom row from left to right with $x_1,\ldots,x_r,y_{r+1},\ldots,y_{r+s}$. For an edge $e$ of $D$ we put $\la e\ra=\la z_1,z_2\ra$, where $z_1$ and $z_2$ are the labels of the endpoints of $e$. Now we define $F(D)=\prod_{e\in D}\la e\ra$.

We have $\GL_m$-equivariant isomorphisms
\begin{align*}
&\End_k(V^{\ot r}\ot(V^*)^{\ot s})\cong\\
&(V^{\ot r}\ot(V^*)^{\ot s})\ot(V^{\ot r}\ot(V^*)^{\ot s})^*\cong\\
&\big(((V^*)^{\ot r}\ot V^{\ot s})\ot(V^{\ot r}\ot(V^*)^{\ot s})\big)^*\cong\\
&\big((V^*)^{\ot (r+s)}\ot V^{\ot(r+s)}\big)^*\, ,
\end{align*}
where the final isomorphism comes from the isomorphism
$$(V^*)^{\ot (r+s)}\ot V^{\ot(r+s)}\stackrel{\sim}{\to}((V^*)^{\ot r}\ot V^{\ot s})\ot(V^{\ot r}\ot(V^*)^{\ot s})$$
that sends $y_1\ot\cdots\ot y_{r+s}\ot x_1\ot\cdots\ot x_{r+s}$ to $$(y_1\ot\cdots\ot y_r\ot x_{r+1}\ot\cdots\ot x_{r+s}\ot x_1\ot\cdots\ot x_r\ot y_{r+1}\ot\cdots\ot y_{r+s})\,.$$
Note that $\big((V^*)^{\ot (r+s)}\ot V^{\ot(r+s)}\big)^*$ can be identified with the vector space of $2(r+s)$-multilinear functions on $\bigoplus^{r+s}V^*\oplus\bigoplus^{r+s}V$. Under the above isomorphisms, $F(D)$ corresponds to an endomorphism $E(D)\in\End_k(V^{\ot r}\ot(V^*)^{\ot s})$ and one can check that $D\mapsto E(D)$ is a representation of $B_{r,s}(m)$. It is clear that the actions of $B_{r,s}(m)$ and $\GL_m$ commute, since the multilinear functions $F(D)$ are $\GL_m$-invariant.

For the proof of Theorem~\ref{thm.doublecentraliser} we need the diagram $b\in B_{r,s}(m)$ which is defined by
\ \\
\begin{equation}
\xy
(-78,3)*{b=},
(-35.88,2.8)*{
\xymatrix @R=14pt @C=14pt @M=-2pt{
{\bullet}{\ar@{-}@/^.7pc/[0,4]}&{\bullet}\ar@{-}[1,0]&\cdots&{\bullet}\ar@{-}[1,0]&{\bullet}&{\bullet}\ar@{-}[1,0]&\cdots&{\bullet}\ar@{-}[1,0]\\
{\bullet}{\ar@{-}@/_.7pc/[0,4]}&{\bullet}&\cdots&{\bullet}&{\bullet}&{\bullet}&\cdots&{\bullet}
}},
(-25,5)*{
\xymatrix @R=40pt @C=14pt @M=-2pt{
{}\ar@{--}[1,0]\\
{}
}}
\endxy\,.
\end{equation}
\vspace{.3cm}

\noindent The main result of the theorem below as the second assertion in (i). If $k$ has characteristic zero this is a trivial consequence of the first assertion, the semisimplicity of the enveloping algebra of $\GL_m$ (by Weyl's theorem), and the fact that a semisimple subalgebra of a matrix algebra equals the centraliser of its centraliser. The theorem was first obtained in \cite[Thm.~6.11]{DiDoSt}.

\begin{thm}\label{thm.doublecentraliser}
The following holds.
\begin{enumerate}[{\rm(i)}]
\item The algebra $\End_{\GL_m}(V^{\ot r}\ot(V^*)^{\ot s})$ coincides with the image of $B_{r,s}(m)$ in $\End_k(V^{\ot r}\ot(V^*)^{\ot s})$ and $\End_{B_{r,s}(m)}(V^{\ot r}\ot(V^*)^{\ot s})$ is the enveloping algebra of $\GL_m$ in $\End_k(V^{\ot r}\ot(V^*)^{\ot s})$.
\item The homomorphism $B_{r,s}(m)\to\End_k(V^{\ot r}\ot(V^*)^{\ot s})$ is injective if and only if $m\ge r+s$.
\end{enumerate}
\end{thm}

\begin{proof}
(i).~The first assertion is proved precisely as in \cite{Br} using the first fundamental theorem (FFT) of invariant theory \cite[Thm.~3.1]{DeCP} for vectors and covectors. Since the group scheme $(\GL_m)_K$ is reduced and $\GL_m$ is dense in $\GL_m(K)$, the $\GL_m$-invariants are the same as the formal invariants, i.e. the invariants of the group scheme $(\GL_m)_k$. The FFT gives us that the space of $2(r+s)$-multilinear $\GL_m$-invariant functions on $\bigoplus^{r+s}V^*\oplus\bigoplus^{r+s}V$ is spanned by the monomials in the $\la x_i,y_j\ra$ with the property that each $x_i$ and each $y_j$ occurs exactly once. These are precisely the functions $F(D)$. Therefore $\End_{\GL_m}(V^{\ot r}\ot(V^*)^{\ot s})$ is spanned by the endomorphisms $E(D)$. One can find similar arguments in \cite{Koike}.

Now we prove the second assertion. Since $B_{r,s}(m)$ is generated by $b$ and the diagrams corresponding to the permutations in $\Sym_r\times\Sym_s$ we have, by Proposition~\ref{prop.envalg}(ii), that $u\in\End_k(V^{\ot r}\ot(V^*)^{\ot s})$ commutes with the action of $B_{r,s}(m)$ if and only if it occurs in the enveloping algebra of $\GL_m\times\GL_m$ and commutes with $b$. It is easy to check that commuting with $b$ amounts to the equations \eqref{eq.envalg}. So the assertion follows from Proposition~\ref{prop.envalgeqs}.\\
(ii).~If $m\ge r+s$, then the second fundamental theorem \cite[Thm.~3.4]{DeCP} gives us that the functions $F(D)$, and therefore the endomorphisms $E(D)$, are linearly independent, since the functions $\la x_i,y_j\ra$, $1\le i,j\le r+s$, are then algebraically independent. If $m<r+s$, then the equation $\det\big(\la x_i,y_j\ra_{1\le i,j\le r+s}\big)=0$ produces a nontrivial linear relation between the $F(D)$ and therefore also one between the $E(D)$.
\end{proof}

\begin{rem}
Using Theorem~\ref{thm.doublecentraliser} one can construct a ``rational Schur functor", compare \cite{DT}, and show, for example, that the walled Brauer algebra and the rational Schur algebra have the same block relation, compare \cite[Thm.~5.5]{DT}.
\end{rem}

\noindent{\it Acknowledgement}. This research was funded by a research grant from The Leverhulme Trust.

\bigskip

{\sc\noindent Department of Mathematics,
University of York, Heslington, York, UK, YO10~5DD.
{\it E-mail address : }{\tt rht502@york.ac.uk}
}

\end{document}